\documentclass{amsart}
\usepackage{amssymb, amsmath, amsfonts, amscd}

\input xypic

\theoremstyle{plain}
\newtheorem{theorem}{Theorem}[section]
\newtheorem{corollary}[theorem]{Corollary}
\newtheorem{lemma}[theorem]{Lemma}
\newtheorem{proposition}[theorem]{Proposition}

\newtheorem{example}[theorem]{Example}

\theoremstyle{definition}
\newtheorem{definition}[theorem]{Definition}

\theoremstyle{remark}

\numberwithin{equation}{theorem}

\newcommand{\m}{\mathfrak{m}}
\newcommand{\F}{\mathcal{F}}

\newcommand{\I}{\mathcal{I}}
\newcommand{\J}{\mathcal{J}} 

\newcommand{\C}{\mathcal{C}}
\newcommand{\Q}{Q}
\renewcommand{\L}{\mathcal{L}}
\newcommand{\E}{\mathcal{E}}

\newcommand{\Hom}{\operatorname{Hom}} 
\renewcommand{\O}{\mathcal{O} }

\renewcommand{\P}{\mathbf{P} }
\renewcommand{\Pr}{\mathcal{J} }

\renewcommand{\equiv}{\backsim }

\newcommand{\SL}{\operatorname{SL}} 
\newcommand{\GL}{\operatorname{GL}}

\renewcommand{\H}{\operatorname{H} }
\newcommand{\sym}{\operatorname{Sym}}

\newcommand{\R}{\operatorname{R} }

\newcommand{\lp}{\mathfrak{p}} 
\renewcommand{\lq}{\mathfrak{q}}
\newcommand{\lm}{\mathfrak{m}}

\newcommand{\Z}{\mathbf{Z} }

\newcommand{\Spec}{\operatorname{Spec} }
\newcommand{\Pic}{\operatorname{Pic} }
\newcommand{\ft}{\frac{1}{t} }

\newcommand{\s}{\mathcal{S}}

\newcommand{\p}{\mathbb{P}} 

\newcommand{\g}{\mathbb{G}}
\renewcommand{\R}{\mathbf{R} }

\newcommand{\tphi}{\tilde{\phi}}
\newcommand{\schur}{\mathbb{S}}
\renewcommand{\k}{\kappa}

\begin{document}

\title{Discriminants of morphisms of sheaves}

\author{Helge Maakestad}

\email{\text{h\_maakestad@hotmail.com} }
\thanks{Partially supported by a research scholarship from www.nav.no}


\date{23.9.2009}

\begin{abstract}  The aim of this paper is to give a unified
  definition of a large class of discriminants arising in algebraic
  geometry using the discriminant of a morphism of 
  locally free sheaves. The discriminant of a morphism of locally free
  sheaves has a geometric definition in terms of
  grassmannian bundles, tautological sequences and projections and is a simultaneous
  generalization of the discriminant of a morphism of schemes, the
  discriminant of a linear system on a smooth projective scheme and
  the classical discriminant of degree $d$ polynomials.
 We study the discriminant of a morphism in various situations: The
 discriminant of a finite morphism of schemes, the discriminant of a
 linear system on the projective line and the discriminant of a linear
 system on a flag variety. The main result of the paper is that the
 discrimiant of any linear system on any flag variety is irreducible.

\end{abstract}

\maketitle

\tableofcontents

\section{Introduction} 

The aim of this paper is to give a unified definition of a large
class of discriminants arising in algebraic geometry.
We define the discriminant of an arbitrary
morphism of locally free finite rank  sheaves. This discriminant is
defined for any morphism 
\[ \phi:u^*\E \rightarrow \F \]
of locally free finite rank $\O_X$-modules on an arbitrary scheme $X$ defined over an arbitrary base
scheme $S$ where $u:X\rightarrow S$ is any  quasi compact morphism of schemes. 
The discriminant of a morphism of locally free  sheaves
has a geometric definition in terms of grassmannian bundles,
tautological sequences and projections.
It is a simultaneous generalization of the discriminant of a morphism of
schemes, the discriminant of a linear system on a smooth projective
scheme and the classical discriminant of degree $d$ polynomials. 
We study this discriminant in the case of linear systems
on projective spaces and flag varieties and use a previous result on the $P$-module
structure of the jet bundle (see \cite{maa1}) to give a candidate for a resolution of
the ideal sheaf $\I$ of the $k$'th discriminant $D^k(\O(d))$ of the
linear system defined by $\O(d)$ on $\p(V^*)$. 
We also study the discriminant of a linear system on any flag variety
$\SL(E)/P$ and prove it is irreducible in general.

In section two
we give the general definition of the discriminant of a morphism of
locally free sheaves. We prove the following general result: Assume we are
given a morphism $\phi:u^*\E\rightarrow \F$ of locally free sheaves on $X$ where $X$ is
irreducible and quasi compact over $S$ and $Coker(\phi^*)$ is locally
free. It follows the discriminant $D^1(\phi)$ is irreducible. 
This is Corollary \ref{irreducible}. 
We also prove it is a simultaneous generalization of the
discriminant of a morphism of schemes and the discriminant of a linear
system on a smooth projective scheme (see Example \ref{morphisms} and \ref{discriminant}).
We prove in Example \ref{classical} that the discriminant $D^1(\O(d))$
on the projective line is the classical discriminant parametrizing
degree d polynomials with multiple roots. 

In section three of the paper we consider the discriminant
$Discr(P(t))$ of a polynomial $P(t)$ in $A[t]$ where $A$ is any
commutative ring. We relate the discriminant to properties of the ring
extension $A\subseteq A[t]/P(t)$ and give precise criteria for the
extension to be etale in the case when $P(t)$ is a monic polynomial.

In section four we study the discriminant $D^l(\O(d))$ of a line
bundle $\O(d)$ on $\p^1_K$ where $K$ is an arbitrary field. Using the
Taylor morphism, jet bundles and projections we prove in Theorem
\ref{main2} that $D^l(\O(d))$ is an irreducible local complete
intersection for all $1\leq l \leq d$.

In section five of the paper we prove some general results on jet
bundles on projective space and higher cohomology groups of exterior
powers of jet bundles on projective space. We give a complete
description (see Theorem \ref{exterior}) of the
$\SL(V)$-module structure of all higher cohomology groups of all
exterior powers of jet bundles and dual jet bundles on projective space.
We also calculate the higher direct images of a class of twisted jet
bundles with an $\SL(V)$-linearization (see Theorem \ref{directimage}).
We also prove that any $G$-module $W$ may be realized as the global
sections of a $G$-linearized locally free $\O_{G/P}$-module $\E(\rho)$ 
(see Proposition \ref{global}). We study the discriminant of a
linear system on projective space and prove existence of a complex of
locally free sheaves which is a candidate for a resolution of the ideal
sheaf of the discriminant (see Example \ref{complex}). 
We finally consider discriminants of linear systems on flag varieties
(see Example \ref{flag}) and prove in Theorem \ref{irr} all such
discriminants are irreducible.

Much research has been devoted to the study of discriminants and
syzygies of discriminants (see \cite{gelfand} and \cite{weyman}). The
novelty of the approach in this paper is the introduction of a
functorial discriminant valid for a map of locally free sheaves relative
to a quasi compact family of schemes. This gives a unified definition of a
large class of discriminants appearing in algebraic geometry. 
All definitions are intrinsic and all discriminants have a canonical
scheme structure.

\section{Discriminants of morphisms of sheaves}

In this section we introduce the discriminant of an arbitrary morphism 
of locally free sheaves on an arbitrary scheme. We prove the
discriminant of a morphism of locally free sheaves is a simultaneous generalization
of the discriminant of a morphism of schemes and the discrimimant of a
linear system on a smooth projective scheme.

Let in the following $u:X\rightarrow S$ be an arbitrary quasi compact morphism of
schemes and let $\E$ be a locally free finite rank $\O_S$-module and $\F$ a
locally free finite rank $\O_X$-module. Assume 
\[ \phi:u^*\E\rightarrow \F \]
is an arbitrary morphism of $\O_X$-modules. 
Let $\g_n(u^*\E^*)$ be the grassmannian bundle of the locally free
sheaf $u^*\E$. The grassmannian has the following properties:
There is a projection morphism
\[ \tilde{\pi}: \g_n(u^*\E^*)\rightarrow X \] 
with the following properties: Let $Y=\g_n(u^*\E^*)$.
There is an isomorphism
\[ Y\cong \g_n(\E^*)\times_S X \]
giving a commutative diagram
\[
\diagram   Y \rto_p \dto^q & X \dto^u \\
          \g_n(\E^*) \rto^{\pi} & S .
\enddiagram \]
Here $p=\tilde{\pi}$ is the projection morphism.

There is on $\g_n(\E^*)$ a \emph{tautological sequence}
\begin{equation} \label{taut}
 0\rightarrow \s \rightarrow \pi^*\E \rightarrow \Q \rightarrow 0 
\end{equation}
of locally free sheaves with $rk(\s)=n$. The locally free sheaf $\s$
is the \emph{tautological subbundle} on $\g_n(\E^*)$.
The sequence \ref{taut} reflects the fact that the
grassmannian $\g_n(\E^*)$ is the scheme representing  the grassmannian
functor $Grass_n(\E^*)$: Via the \emph{Yoneda Lemma} it follows the
grassmannian functor $Grass_n(\E^*)$ is represented by a scheme
$\g_n(\E^*)$ and a \emph{universal object}. The universal object is
given by the subbundle
\[ 0\rightarrow \s \rightarrow \pi^*\E .\]

\begin{example} The tautological line bundle.
\end{example}
In the case when $n=1$ it follows $\s=\O(-1)$ and we get the
sequence of the tautological sub-bundle
\[ 0\rightarrow \O(-1)\rightarrow \pi^*\E\]
on $\p(\E^*)$.

\begin{proposition} \label{quasicompact} Let $u:X\rightarrow S$ be a
  quasi compact morphism
  of schemes and let $\E$ be a locally free $\O_S$-module of rank
  $m$. It follows
\[ \pi:\g_n(\E^*) \rightarrow S \]
and
\[ q:\g_n(u^*\E^*)\rightarrow \g_n(\E^*) \]
are  quasi compact morphisms.
\end{proposition}
\begin{proof} Since $u$ is quasi compact it follows for any open
  affine subscheme $V=\Spec(A)\subseteq S$ the inverse image
  $U=u^{-1}(V)\subseteq X$ is a finite union of open affine schemes:
\[U=\Spec(B_1)\cup \cdots \cup \Spec(B_k) .\]
This is a general fact: A scheme over an affine scheme is quasi
compact if and only if it is a finite union of open affine sub
schemes.
Let $\pi_V:V\rightarrow \Spec(\Z)$ and $\pi_U:U\rightarrow \Spec(\Z)$
be the structure morphisms. Assume $E|_V=A\{e_1,..,e_m\}=\pi_V^*W$
where
\[ W=\Z\{e_1,..,e_m\}.\]
Let
\[ W^*=\Z\{x_1,..,x_m\}.\]
It follows
\[ \pi^{-1}(V)=\g_n(\E^*|_V)\cong \g_n(\pi_V^*W^*)\cong
\g_n(W^*)\times_\Z \Spec(A).\]
Since $\g_n(W^*)=\cup_{i=1}^l \Spec(A_i)$ 
is a finite union of affine open schemes $\Spec(A_i)$ it follows
\[ \pi^{-1}(V)=\g_n(W^*)\times_\Z \Spec(A)=\cup_{i=1}^l
\Spec(A_i)\times_\Z \Spec(A)\cong \]
\[  \cup_{i=1}^l \Spec(A_i \otimes_\Z A)   \]
is a finite union of affine open schemes. It follows $\pi$ is a quasi
compact morphism of schemes.

Pick the open set
\[ \Spec(A_i\otimes_\Z A)\subseteq \pi^{-1}(V)=\g_n(W^*)\times_\Z
\Spec(A).\]
It follows
\[ q^{-1}(\Spec(A_i\otimes_\Z A))=\Spec(A_i)\times_V u^{-1}(V)= \]
\[ \cup_{j=1}^k \Spec(A_i)\times_V \Spec(B_j)=\cup_{j=1}^k
\Spec(A_i\otimes_\Z B_j)  .\]
Hence
\[ q^{-1}(\Spec(A_i\otimes_\Z A)) \]
is a finite union of open affine schemes.
The open sets $\Spec(A_i\otimes_\Z A)$ cover $\g_n(W^*)\times_\Z
\Spec(A)$ hence it
follows $q$ is a quasi compact morphism, and the Proposition is
proved.
\end{proof}

Note: For any quasi compact morphism $u:X\rightarrow S$ of schemes and
any closed subscheme $Z\subseteq X$ we get an induced morphism
$v:Z\rightarrow S$ which is quasi compact. The map of structure
sheaves
\[ v^\#:\O_S\rightarrow v_*\O_Z \] gives rise to an ideal sheaf
\[ \I=ker(v^\#)\subseteq \O_S .\]
The ideal sheaf $\I\subseteq \O_S$ corresponds to a subscheme
$u(Z)\subseteq S$: the \emph{schematic image of $Z$ via $u$}. Hence if
$v:Z\rightarrow S$ is a closed morphism it follows we get a canonical
structure of closed subscheme on the schematic image $u(Z)\subseteq
S$. This structure is not neccessarily the reduced induced structure
on
the topological space $u(Z)$ viewed as a closed subspace of $S$.

The sequence \ref{taut} has the following property: Assume $s\in \pi^{-1}(z)$ is a
$\kappa(z)$-rational point. If we take the fiber of \ref{taut} at $s$
we get an exact sequence of $\kappa(z)$-vector spaces
\[ 0\rightarrow \s(s) \rightarrow \pi^*\E(s) \rightarrow
\Q(s)\rightarrow 0 \]
and $\pi^*\E(s)=\E(\pi(s))\otimes_{\kappa(\pi(s))}\kappa(s)\cong \E(z)$.
We get a canonical $n$-dimensional $\kappa(z)$ sub vector space
\[ S(s) \subseteq \E(z) \]
for each $s\in \pi^{-1}(z)(\kappa(z))$.
By functoriality the following holds:
\[ \pi^{-1}(z)\cong \g_n(\E(z)^*) \]
hence the fiber $\pi^{-1}(z)$ is canonically isomorphic to the
grassmannian parametrizing $n$-planes in the $\kappa(z)$-vector space
$\E(z)$. The tautological sequence \ref{taut} gives for each
$\kappa(z)$-rational point in $\g_n(\E(z)^*)$ its corresponding
$n$-plane in $\E(z)$. We get a one to one correspondence 
\[\{ \k(z)-\text{rational points }s\in\g_n(\E(z)^*)\} \cong \{ n-\text{planes }W\subseteq \E(z)\} \]
given by
\[ s \in \g_n(\E(z)^*)\equiv \s(s)\subseteq \E(z).\]

We will use the tautological sequence to define the discriminant of a
morphism of locally free sheaves. Let $Y=\g_n(\E^*)\times_S X$.
 Let $\s_Y=q^*\s$, $\E_Y=q^*\pi^*\E=p^*u^*\E$ and $\F_Y=p^*\F$. We get
a morphism
\[ \s_Y \rightarrow \E_Y \rightarrow^{p^*\phi} \F_Y . \]
Let $\tilde{\phi}$ be the composed morphism
\[ \tilde{\phi}:\s_Y \rightarrow \F_Y.\]
Let $Z(\tilde{\phi})$ be the zero scheme of the morphism $\tilde{\phi}$. If
$U=\Spec(A)$ is a trivialization of $\s_Y$ and $\F_Y$ and
$\tilde{\phi}=(a_{ij})$ with $a_{ij}\in A$ it follows the ideal of
$Z(\tilde{\phi})$ is generated by $a_{ij}$ on the open set $U$.

Since $\pi$ and $q$ are quasi compact morphisms we may define the following:

\begin{definition}  Let $I^n(\phi)=Z(\tilde{\phi})$ be the \emph{n-incidence scheme} of $\phi$.
The scheme $D^n(\phi)=q(I(\phi))$  is the \emph{n-discriminant} of $\phi$. 
The scheme $Discr^n(\phi)=\pi(D^n(\phi))$ is the \emph{direct image n-discriminant} of
  the morphism $\phi$.
\end{definition}

We get a diagram
\[
\diagram I^n(\phi) \rto^i \ddto^{\tilde{q}} & \g_n(\E^*)\times_S X \rto_p \dto^q & X \dto^u \\
           &    \g_n(\E^*) \rto^{\pi} & S \\
        D^n(\phi) \rrto \urto^j & & Discr^n(\phi) \uto^k
\enddiagram \]
where $i,j$ and $k$ are inclusions of schemes.

Assume $\psi:\E\rightarrow \F$ is a map of locally free $\O_X$-modules.
\begin{lemma} \label{zero}The following holds: $x\in Z(\psi)$ if and
  only if $\psi(x)=0$.
\end{lemma}
\begin{proof} Let $U=\Spec(A)\subseteq X$ be an open subset where $\E$
  and $\F$ trivialize. It follows $\E$ is the sheafification of $A^m$
  and $\F$ the sheafification of $A^n$ for some integers $m,n\geq 1$. 
Let $x\in \Spec(A)$ correspond to a prime ideal $\lp_x\subseteq A$.
It follows $x\in Z(\psi)$ if and only of $a_{ij}\in \lp_x$ for all
$i,j$ where
$(a_{ij})=\psi|_U$ and $a_{ij}\in A$ are the coefficients of $\psi$
over $U$. It follows $x\in Z(\psi)$ if and only if the fiber map
$\psi(x)$ is zero and the Lemma is proved.
\end{proof}

Let $\phi^*:\F^*\rightarrow u^*\E^*$ be the dual of $\phi$ and
consider the exact sequence
\begin{equation} \label{cok}
\F^*\rightarrow^{\phi^*} u^*\E^* \rightarrow
Coker(\phi^*)\rightarrow 0 
\end{equation}
of coherent $\O_X$-modules. Assume $Coker(\phi^*)$ is locally free and
let $n=1$. We get a closed immersion
\[ \p(Coker(\phi^*)) \subseteq \p(u^*\E)\cong \p(\E^*)\times_S X\]
of schemes.

\begin{theorem} \label{equality} There is an equality $I^1(\phi)=\p(Coker(\phi^*))$ as
  subschemes of $\p(u^*\E^*)$.
\end{theorem}
\begin{proof} Let $Y=\p(u^*\E^*)=\p(\E^*)\times_S X$ and let $\I,\J\subseteq \O_Y$ be the
  ideal sheaves of $I^1(\phi)$ and $\p(u^*\E)$. We want to prove there
  is an equality $\I=\J$ of ideal sheaves. 
Consider the diagram
\[
\diagram Y \rto^p \dto^q    &   X \dto^u \\
               \p(\E^*) \rto^{\pi} &   S
\enddiagram 
\]
Let
\[ \alpha:\O(-1) \rightarrow \pi^*\E \]
be the tautological sub-bundle on $\p(\E^*)$. Pull this and the
morphism
\[ \phi:u^*\E \rightarrow \F \]
back to $Y$ to get the morphism
\[ \tphi: \O(-1)_Y \rightarrow \E_Y \rightarrow \F_Y .\]
By definition $Z(\tphi)=I^1(\phi)$. We want to show
$Z(\tphi)=\p(Coker(\phi^*))$ is an equality of schemes. We prove there
is an equality of ideal sheaves. Assume $V=\Spec(A)\subseteq S$ is an
affine open subscheme where $\E$ trivialize. Let
$\pi_V:\Spec(A)\rightarrow \Spec(\Z)$ be the structure morphism and
let
\[ W=\Z\{e_0,..,e_l\} .\]
It follows 
\[ E|_V=\pi_V^*W=A\otimes_\Z \Z\{e_0,..,e_l\} .\]
Let $U=\Spec(B)\subseteq u^{-1}(V)$ be an open set where $\F$
trivialize and let $\pi_U:\Spec(B)\rightarrow \Spec(\Z)$ be the
structure morphism. Let $Z=\Z\{f_0,..,f_m\}$. It follows
\[ \F|_U=\pi_U^*Z=B\otimes_\Z \Z\{f_0,..,f_m\} .\]
Let $x_i=e_i^*$ and $y_j=f_j^*$.
Pull $\E$ back to $U$ to get
\[ u^*\E|_U=B\otimes_\Z W=B\otimes_Z
\Z\{e_0,..,e_l\}=B\{e_0,..,e_l\}.\]
Restrict the morphism $\phi$ to $U$ to get
\[ \phi|_U:B\{e_0,..,e_l\}\rightarrow B\{f_0,..,f_m\} \]
with $\phi|_U=(b_{ij})$ with $b_{ij}\in B$. Consider the morphism
\[ p:\p(u^*\E^*)\rightarrow X .\]
It follows
\[ p^{-1}(U)=\p(u^*\E|_U)=\p(W^*)\times_\Z \Spec(B).\]
Hence $p^{-1}(U)$ may be covered by open affine schemes on the form
\[U_i= D(x_i)\times_\Z\Spec(B)=\Spec(B[\frac{x_0}{x_i},..,\frac{x_l}{x_0}]) \]
for $i=0,..,l$. Let $t_j=x_j/x_i$ and $t_i=1$.
On $U_i$ the map $\alpha|_{U_i}$ looks as follows:
\[\alpha|_{U_i}:B[\frac{x_0}{x_i},..,\frac{x_l}{x_i}]\frac{1}{x_i}\rightarrow 
B[\frac{x_0}{x_i},..,\frac{x_l}{x_i}]\otimes_\Z\{e_0,..,e_l\} \]
with
\[ \alpha|_{U_i}(1/x_i)=t_0\otimes e_0+t_1\otimes e_1+\cdots +
1\otimes e_i +\cdots +t_l\otimes e_l=\]
\[ [t_0,t_1,..,t_{i-1},1,t_{i+1},..,t_l].\]
It follows the composed morphism $\tphi|_{U_i}\circ \alpha$ has
coefficients on the form
\[ \{ c_{k,i}=b_{k,0}t_0+b_{k,1}t_1+\cdots b_{k,i}+\cdots
+b_{k,l}t_l\}_{k=0}^m .\]
where $b_{i,j}\in B$ are the coefficients of $\tphi|_{U_i}$.
Hence the ideal sheaf $\I$ is generated by the elements $c_{k,i}$ on
the open set $U_i=D(x_i)\times U\subseteq \p(u^*\E^*)$. Consider the
exact sequence
\[ \F^*\rightarrow^{\phi^*}u^*\E^* \rightarrow
Coker(\phi^*)\rightarrow 0.\]
We want to calculate generators for the image $Im(\phi^*)\subseteq
u^*\E^*$ on the open set $U_i$. The matrix of $\phi^*|_{U_i}$ is the
transpose of the matrix $\phi|_{U_i}$ and one checks that on $U_i$ the
following holds:
\[ \phi^*|_{U_i}(y_k)=b_{k,0}t_0+b_{k,1}t_1+\cdots b_{k,i}+\cdots
+b_{k,l}t_l  =c_{k,i}. \]
It follows $Im(\phi^*|_{U_i})$ is generated by the elements $c_{k,i}$
hence the ideal sheaf $\J$ of $\p(Coker(\phi^*))$ is on $U_i$
generated by $c_{k,i}$. It follows $\I=\J$ and the claim of the
Theorem follows.
\end{proof}

\begin{corollary} \label{irreducible} If $X$ is irreducible and $Coker(\phi^*)$ is locally
  free it follows $D^1(\phi)$ is irreducible.
\end{corollary}
\begin{proof} Since $X$ is irreducible and $I^1(\phi)=\p(Coker(\phi^*))$
  is a projective bundle on $X$ it follows $I^1(\phi)$ is
  irreducible. It follows $D^1(\phi)=q(I^1(\phi))$ is irreducible and
  the claim of the Corollary follows.
\end{proof}

We give an interpretation of $I^n(\phi), D^n(\phi)$ and
$Discr^n(\phi)$ in terms of points. Assume $z\in S$ is a point with
residue field $\k(z)$. Assume $y=(s,x)\in Y$ with $u(x)=\pi(s)=z$. It follows $\k(y)=\k(s)=\k(x)=\k(z)$.
By Lemma \ref{zero} it follows $y\in I^n(\phi)$ if and only if
$\phi_Y(y)=0$.
We get
\[ \s_Y(y)=\s(s)\otimes_{\k(s)}\k(y)\cong \s(s) \]
\[ \E_Y(y)=\E(z)\otimes_{\k(z)}\k(y)\cong \E(z) \]
and
\[ \F_Y(y)=\F(x)\otimes_{\k(x)}\k(y)\cong \F(x). \]
It follows the composed map
\[ \tilde{\phi}(y):\s(s)\rightarrow \E(z)\rightarrow^{\phi(x)}\F(x) \]
is the zero map. Hence a $\k(z)$-rational point $y=(s,x)$ is in
$I^n(\phi)$ if and only if the canonical inclusion
\[ \s(s)\subseteq \E(z) \]
induces an inclusion 
\[ \s(s)\subseteq Ker(\phi(x)).\]

\begin{lemma} If $y=(s,x)\in I^n(\phi)(\k(z))$ it follows
  $dim_{\k(z)}Ker(\phi(x))\geq n$.
\end{lemma}
\begin{proof} By the discussion above there is an inclusion
\[ \s(s)\subseteq Ker(\phi(x)) \]
and since $dim_{\k(z)}\s(s)= n$ the claim of the lemma is proved.
\end{proof}

Hence $s\in D^n(\phi)$ with $\pi(s)=z$ if and only if there is a point
$x\in X$ with $u(x)=\pi(s)$ such that the canonical inclusion
\[ \s(s)\subseteq \E(z) \]
induce an inclusion 
\[ \s(s)\subseteq Ker(\phi(x)).\]

\begin{example} On the case $n=1$:
\end{example}
Consider the tautological subbundle on $\p(\E^*)=\g_1(\E^*)$
\[ 0\rightarrow \O(-1)\rightarrow \pi^*\E .\]
Pull this sequence back to $Y=\p(\E^*)\times X$ to get a sequence
\[ \O(-1)_Y\rightarrow \E_Y.\]
Similarly pull back the sequence $\phi:u^*\E\rightarrow \F$ to $Y$ to
get the sequence
\[ \tilde{\phi}: \O(-1)_Y\rightarrow \E_Y\rightarrow \F_Y.\]
Since
\[ \tilde{\phi}\in \Hom(\O(-1)_Y,\F_Y)\cong \Hom(\O,\O(1)_Y\otimes
\F_Y)=\H^0(Y,\O(1)_Y\otimes \F_Y) \]
It follows we have described the $1$-incidence scheme $I^1(\phi_Y)$ as the
zero locus of $\tilde{\phi}$ viewed as a global section of the locally free sheaf
$\O(1)_Y\otimes \F_Y$.

\begin{lemma} There is an exact sequence on $Y$
\[ \O(-1)_Y\otimes \F_Y^*\rightarrow \O_Y \rightarrow
\O_{Z(\tilde{\phi})}\rightarrow 0 .\]
\end{lemma}
\begin{proof} The proof is left to the reader as an exercise.
\end{proof}

Consider for $j\geq 1$ the $j$'th exterior product
\[ \wedge^j(\O(-1)_Y\otimes \F_Y^*)\cong \O(-j)_Y\otimes \wedge^j \F^*_Y.\]
The ideal sheaf $\I$ of $Z(\tilde{\phi})$ is locally generated by a
regular sequence hence from \cite{altman} we get a Koszul-resolution
\[0\rightarrow \O(-r)_Y\otimes \wedge^r \F^*_Y\rightarrow \O(-r+1)_Y\otimes
\wedge^{r-1}\F^*_Y\rightarrow \cdots \]
\[ \cdots \rightarrow \O(-j)_Y\otimes \wedge^j\F^*_Y\rightarrow \cdots \]
\[ \cdots \O(-2)_Y\otimes \wedge^2\F^*_Y\rightarrow \O(-1)_Y\otimes
\F^*_Y \rightarrow \O_Y \rightarrow \O_{Z(\tilde{\phi})}\rightarrow 0\]
of the structure sheaf $\O_{Z(\tilde{\phi})}$ of the incidence scheme $Z(\tilde{\phi})$.
Here $r=rk(\F^*)$.

\begin{proposition} There is an isomorphism of $\O_{\p(\E^*)}$-modules
\[ \R^iq_*(\O(-j)_Y\otimes \wedge^j\F^*_Y)\cong \O(-j)\otimes \pi^*\R^iu_*(\wedge^j\F^*).\]
\end{proposition}
\begin{proof} By the projection formula and base change we get
\[\R^iq_*(\O(-j)_Y\otimes \wedge^j\F^*_Y)\cong
\R^iq_*(q^*\O(-j)\otimes \wedge^jq^*\F^*)\cong \]
\[ \O(-j)\otimes \R^iq_*(p^*\wedge^j\F^*)\cong\O(-j)\otimes\pi^*
\R^iu_*(\wedge^j\F^*)\]
and the claim of the Proposition is proved.
\end{proof}

We get a double complex of $\O_{\p(\E^*)}$-modules on $\p(\E^*)$
defined by
\[ \C(\phi)^{i,j}=\O(-j)\otimes \pi^*\R^iu_*(\wedge^j\F^*) .\]
A natural question to ask is if the total complex
\begin{equation} \label{total} Tot(\C(\phi)^{i,j})_n=\oplus_{i+j=n}\O(-j)\otimes
  \pi^*\R^iu_*(\wedge^j\F^*) 
\end{equation}
may be used to construct a resolution of the ideal sheaf of $D(\phi)\subseteq \p(\E^*)$.
\begin{definition} Let the total complex $Tot(\C(\Phi)^{i,j})_n$ be
  the \emph{discriminant complex} of $\phi$.
\end{definition}

\begin{example} Discriminants of morphisms of schemes.
\label{morphisms}
\end{example}
Assume $f:U\rightarrow V$ is a quasi compact map of smooth schemes of finite type
over a field $F$. We get a map of sheaves of differentials
\begin{equation} \label{cotangent} df: f^*\Omega^1_{V}\rightarrow  \Omega^1_{U} .
\end{equation}
Since $U,V$ are smooth over $F$ it follows $\Omega^1_U, \Omega^1_V$
are locally free sheaves of finite rank. 
We use the cotangent sequence \ref{cotangent} to give a set theoretic
definition of the
discriminant of the morphism $f$. 
\begin{definition} \label{settheoretical}
Let $Discr^1(f)\subseteq V$ be the set of points $s\in V$ such that
there is an $x\in U$ with $f(x)=s$ and $Ker(df(x))\neq 0$. We say
$Discr^1(f)$ is the \emph{set theoretic discriminant of the morphism $f$}.
\end{definition}
Let $I^n(df)$ be the $n$-incidence scheme of $df$. We get a diagram of maps of schemes
where by Proposition \ref{quasicompact} $q$ and $\pi$ are quasi compact morphisms: 
\[
\diagram I^n(df) \rto^i \ddto^{\tilde{q}} & \g_n((\Omega^1_V)^*)\times U \rto_p \dto^q & U \dto^f \\
           & \g_n((\Omega^1_V)^*) \rto^\pi & V \\
        D^n(df) \urto \rrto & & Discr^n(df) \uto 
\enddiagram.
\]
\begin{definition} The scheme $D^n(df)$ is the \emph{$n$-discriminant} of the
  morphism $f$. The scheme $Discr^n(df)$ is the \emph{direct image
    $n$-discriminant} of the morphism $f$.
\end{definition}

Pick $y=(d,s)\in I^n(df)$. Assume $\pi(d)=f(s)=z$ and
$\k(y)=\k(d)=\k(s)=\k(z)$. Since $d\in
\pi^{-1}(z)=\g_n(\Omega^1_V(z)^*)$ the tautological sequence on the
grassmannian gives a canonical $n$-dimensional vector subspace
\[ \s(d)\subseteq \Omega^1_V(z).\]
We get a composed map
\[ \s(d)\rightarrow \Omega^1_V(z)\rightarrow^{df(s)} \Omega^1_U(s) \]
and since the composed map is the zero map it follows
\[ \s(d)\subseteq Ker(df(s)). \] 
It follows $dim_{\k(z)}Ker(df(s))\geq n$ hence $df(s)$ is not
injective at $s$.

It follows the underlying set of points of $Discr^n(df)$ is the set of
$s=f(x)\in V$ with $x\in U$ and $dim(Ker(df(x))\geq n$.
We see the discriminant of a morphism of locally free sheaves
generalize the set theoretic discriminant of a morphism of schemes in the sense that
the underlying set of points of $D^1(df)$ is a lifting of $Discr^1(f)\subseteq V$ to the projectivization
of the cotangent bundle $\Omega^1_V$. We get a sequence of
subschemes
\[ \cdots \subseteq Discr^n(df)\subseteq \cdots \subseteq Discr^1(df)\subseteq
V .\]
In this case the total complex \ref{total}
becomes
\[ Tot(\C(f)^{i,j})_n=\oplus_{i+j=n}\O(-j)\otimes
\pi^*\R^if_*(\wedge^j (\Omega_U^1)^*).\]

\begin{example} \label{discriminant} Discriminants of linear systems on projective schemes.\end{example}
Let $X\subseteq \p^d_F$ be a smooth projective scheme over a field $F$
and let $\L\in \Pic(X)$. Let $u:X\rightarrow \Spec(F)$ be the
structure morphism and let $W=\H^0(X,\L)$. There is a morphism of locally
free sheaves
\[ 
u^*W\rightarrow^{T^k}\Pr^k_X(\L). \]
Here $T^k$ is the $k$'th Taylor map and $\Pr^k_X(\L)$ is the $k$'th
jet bundle of $\L$. We get a diagram of maps of schemes
\[
\diagram I^1(T^k)\rto^i \dto & \p(W^*)\times X \rto_p
\dto^q & X \dto^u \\
          D^1(T^k) \rto^j & \p(W^*) \rto^\pi & \Spec(F)
\enddiagram
\]
\begin{definition}\label{linear}
The scheme $D^k(\L)=D^1(T^k)$ is the \emph{$k$'th discriminant} of the linear system
defined by $\L$.
\end{definition}

We see the discriminant of a map of locally free sheaves
generalize the discriminant of a linear system on a smooth projective
scheme. The classical discriminant $D^k(\L)$ equals the $1$-discriminant $D^1(T^k)$
of the $k$'th Taylor morphism.

In this case the total complex \ref{total} becomes
\[ Tot(\C(T^k)^{i,j})_n=\oplus_{i+j=n}\O(-j)\otimes \pi^*
\H^i(X,\wedge^j \Pr^k(\L)^*) .\]

To study the total complex for discriminants of linear systems we need
information on the higher cohomology of exterior powers of duals of jet bundles.
In the next section we will study higher cohomology groups of exterior
powers of  $\SL(V)$-linearized jet bundles and the total complex in the situation
where $X=\p(V^*)$ and $\L=\O(d)$ for $d\geq 1$.

Note: The definition given in \ref{linear} was communicated to the
author by D. Laksov.

\begin{example} Invertible sheaves on projective space.\end{example}

We interpret the rational points of $D^k(\O(d))$ on $\p(V^*)$ where
$V$ is an $N+1$-dimensional vector space over any field $F$. Let $W=\H^0(\p(V^*),\O(d))$
Consider
the diagram
\[
\diagram   \p(W^*)\times \p(V^*) \rto_p \dto^q & \p(V^*) \dto^\pi \\
           \p(W^*) \rto^\pi &  \Spec(F)
\enddiagram
\]
Let $Y=\p(W^*)\times \p(V^*)$ and consider the
following sequence of locally free $\O_Y$-modules
\[ \phi_Y:\O(-1)_Y\rightarrow W \otimes \O_Y
\rightarrow \Pr^k(\O(d))_Y .\]
It follows the $F$-rational points of the incidence scheme $Z(\phi_Y)$
have the following interpretation:
Pick an $F$-rational point $x=(s,y)\in Y$ $\kappa(x)=\kappa(s,y)=F$. It
follows there is an equality of residue fields $\kappa(s)=\kappa(y)=\kappa(s,y)=F$.
The point $x$ is by Lemma \ref{zero} in $Z(\phi_Y)$ if and only if $\phi_Y(x)=0$. We
interpret this equation in terms of fibers: We get
\[ \O(-1)_Y(s,y)\rightarrow W \otimes \O_Y(s,y)\rightarrow
\Pr^k(\O(d))_Y(s,y) \]
which becomes
\[ \O(-1)(s)\otimes_{\kappa(s)}\kappa(s,y)\rightarrow
W  \rightarrow
\Pr^k(\O(d))(y)\otimes_{\kappa(y)}\kappa(s,y) \]
which becomes
\[ \O(-1)(s)\rightarrow W
\rightarrow^{T^k(y)}\Pr^k(\O(d))(y) .\]

Write $\tilde{s}=\O(-1)(s)\subseteq W$.
It follows 
\[ T^k(y)(\tilde{s})=0\text{ in } \Pr^k(\O(d))(y).\]
Hence  an $F$-rational $(s,y)\in \p(W^*)\times \p(V^*)$ is in
$Z(\phi_Y)$ if and only if its corresponding section
$\tilde{s}\subseteq W$ satisfies
\[ T^k(y)(\tilde{s})=0 .\]
Hence the points $s\in D^k(\O(d))(F)$ are described in terms
of the taylor map $T^k$ at some point $y\in \p(V^*)(F)$. In local
coordinates the Taylor map $T^k(y)$ formally taylor
expands a global section $\tilde{s}\in W$. 

\begin{example} \label{classical} Linear systems on the projective line.
\end{example}

Assume $\p(V^*)=\p^1$ is the projective line where $V=F\{e_0,e_1\}$
and $V^*=F\{x_0,x_1\}$. Let $W=\H^0(\p^1,\O(d))$.
The Taylor map
\[ T^1:W\otimes \O_{\p^1}\rightarrow \Pr^1(\O(d)) \]
is defined as follows: Let $s_i=x_0^{d-i}x_1$ for $i=0,..,d$ be the
global sections of $\O(d)$. Let $y_i=s_i^*, t=\frac{x_1}{x_0}$ and
$\frac{1}{t}=\frac{x_0}{x_1}$. 
Let $U_{ij}=D(y_i)\times
D(x_j)$ be an open cover of $\p(W^*)\times \p^1$. Let
$u_j=\frac{y_j}{y_i}$.
On $D(x_0)$ we get the following 
\[ T^1:F[t]\{s_i\}\rightarrow F[t]\{1\otimes x_0^d,dt\otimes x_0^d\} \]
with
\[ T^1(s_i)=T^1(x_0^{d-i}x_1^i)=t^i\otimes x_0^d+it^{i-1}dt\otimes
x_0^d.\]
On $D(x_1)$ we get
\[ T^1:F[1/t]\{s_i\}\rightarrow F[1/t]\{1\otimes x_1^d,d(1/t)\otimes x_1^d\} \]
with
\[ T^1(s_i)=T^1(x_0^{d-i}x_1^i)=(1/t)^i\otimes
x_1^d+i(1/t)^{i-1}dt\otimes x_1^d.\]
Consider the following map on $Y=\p(W^*)\times \p^1$
\begin{equation}\label{map1}
\O(-1)_Y\rightarrow^{\alpha} W_Y \rightarrow^{T^1_Y} \Pr^1(\O(d))_Y 
\end{equation}
and restrict to $U_{i0}=D(y_i)\times D(x_0)$.
We get
\[ \alpha|_{U_{i0}}:F[u^i_j,t]\frac{1}{y_i}\rightarrow
F[u^i_j,t]\otimes_F F\{s_0,..,s_d\} \rightarrow^{T^1} F[u^i_j,t]\{1\otimes
x_0^d,dt\otimes x_0^d\} \]
given by
\[ T^1(\alpha(1/y_i))=T^1(u^i_0\otimes s_0 +\cdots 1\otimes s_i
+\cdots u^i_d\otimes s_d)= \]
\[ u^i_0\otimes x_0^d+u^i_1(t+dt)\otimes x_0^d+\cdots +(t+dt)^i\otimes
x_0^d+\cdots +u^i_d(t+dt)^d\otimes x_0^d=\]
\[ f_i(t)\otimes x_0^d+f_i'(t)dt\otimes x_0^d \]
where
\[f_i(t)=u^i_0+u^i_1t+\cdots +u^i_dt^d.\]
Let 
\[ a(t)=y_0+y_1t+\cdots +y_dt^d \]
it follows $y_if_i(t)=a(t)$.
Restrict the map \ref{map1} to $U_{i1}=D(y_i)\times D(x_1)$. We get
\[\alpha|_{U_{i1}}:F[u_j^i,1/t]\frac{1}{y_i}\rightarrow
F[u^i_j,1/t]\otimes_F F\{s_i\}\rightarrow^{T^1}F[u^i_j,1/t]\{1\otimes
x_1^d,d(1/t)\otimes x_1^d\} \]
given by
\[T^1(\alpha(1/y_i))=T^1(u^i_0\otimes s_0 +\cdots 1\otimes s_i
+\cdots u^i_d\otimes s_d) =\]
\[u^i_0(\ft)^d\otimes x_1^d+(d-1)(\ft)^{d-1}d(\ft)\otimes
x_1^d)+\cdots +\]
\[ ((\ft)^{d-i}\otimes
x_1^d+(d-i)(\ft)^{d-i-1}d(\ft)\otimes x_1^d)+\cdots + u^i_d\otimes
x_1^d=\]
\[g_i(\ft)\otimes x_1^d+g_i'(\ft)d(\ft)\otimes x_1^d \]
where
\[ g_i(\ft)=u^i_0(\ft)^d+\cdots +(\ft)^{d-i}+\cdots + u^i_d.\]
Let 
\[b(t)=y_0(\ft)^d+\cdots +y_i(\ft)^{d-i}+\cdots +y_d.\]
It follows $b(\ft)=y_ig_i(\ft)$.
Let $U_{ij}=D(y_i)\times D(x_j)$ and let $Y_{ij}=U_{ij}\cap
I^1(\O(d))$. It follows $q(Y_{ij})=Z_{ij} \subseteq D(y_i)$. Let $j=0$ and let
$I_{i0}=(f_i(t),f'_i(t))$. 
We get a map
\[ q:Y_{i0}\rightarrow D(y_i) \]
which gives a map
\[ q^{\#}:\O_{D(y_i)}\rightarrow q_*\O_{Y_{i0}} \]
given by
\[ q^{\#}:F[u^i_j]\rightarrow F[u^i_j,t]/I_{i0} .\]
It follows
\[  \I_{Z_{i0}}=ker(q^{\#})=Res(f_i(t),
f_i'(t))=Res(a(t),a'(t))|_{D(y_i)}.\]
Let $j=1$ and let $I_{i1}=(g_i(\ft),g_i'(\ft))$. We get a map
\[ q:Y_{i1}\rightarrow D(y_i) \]
which gives a map
\[ q^{\#}:\O_{D(y_i)}\rightarrow q_*\O_{Y_{i1}} \]
given by
\[ q^{\#}:F[u^i_j]\rightarrow F[u^i_j,\ft]/I_{i1}.\]
It follows
\[ \I_{Z_{i1}}=ker(q^{\#})=Res(b_i(\ft),b_i'(\ft))|_{D(y_i)}.\]
Let $\I=(Res(a(t),a'(t))$ and $\J=Res(b(\ft),b'(\ft))$. It follows
\[
\I_{D^1(\O(d))}|_{Z_{i0}}=\I_{Z_{i0}}=Res(a(t),a'(t))|_{D(y_i)}=\I|_{D(y_i)}
\]
and
\[ \I_{D^1(\O(d))}|_{Z_{i1}}=\I_{Z_{i1}}=Res(b(\ft),b'(\ft))|_{D(y_i)}=\J|_{D(y_i)}
.\]
We get an equality of ideal sheaves in $\O_{\p(W^*)}$
\[ \I_{D^1(\O(d))}=\I=\J .\]
It follows the ideal sheaf $\I_{D^1(\O(d))}$ is generated by the
irreducible polynomial
\[ Res(a(t),a'(t))=Res(b(\ft),b'(\ft)) \]
hence $D^1(\O(d))$ is the \emph{discriminant
scheme of degree $d$ polynomials} parametrizing degree $d$ polynomials 
\[ a(t)=y_0+y_1t+\cdots +y_dt^d \]
in the variable $t$ with multiple roots. 
It follows $D^1(\O(d))$ is a determinantal scheme.
We get a filtration of closed subschemes
\[ D^d(\O(d)) \subseteq   \cdots   D^i(\O(d))  \cdots \subseteq D^1(\O(d))       \subseteq \p(W^*). \]
One may ask if $D^i(\O(d))$ is a determinantal scheme for
$1<i\leq d$. There is work in progress on this problem: One wants to check
if $D^k(\L(\lambda))$ is a determinantal scheme where $\L(\lambda)\in
\Pic^G(G/P)$. Here $G$ is a semi simple linear algebraic
group and $P\subseteq G$ a parabolic sub group.

Note: For determinantal schemes much is known about their syzygies
(see \cite{lascoux}).

\section{Discriminants and standard etale morphisms}

In this section we study the discriminant of a polynomial and
its relationship with finite and standard etale morphisms. We relate the discriminant of
a monic polynomial $P(t)$ in $A[t]$ to properties of the integral ring extension
$A\subseteq A[t]/P(t)$. We give an explicit proof of the fact that 
\[ \Spec(A[t]/P(t))\rightarrow \Spec(A) \]
is generically etale using the discriminant $Discr(P(t))$.
We prove in Theorem \ref{standard} that any etale morphism of schemes
is locally on the form
\[ \Spec(A[t]/P(t))\rightarrow \Spec(A) \]
where $P(t)\in A[t]$ is a polynomial with $Discr(P(t))$ a unit in $A$.
We also prove a general result (see Theorem \ref{main1}) on properties
of the ring extension $A\subseteq A[t]/P(t))$ when $P(t)$ is an arbitrary (not necessarily monic) polynomial.

Let $X=\Spec(B)$ and $Y=\Spec(A)$ where $A,B$ are commutative unital rings.
Let $f:\Spec(B)\rightarrow \Spec(A)$ be a finite map of affine
schemes. Hence $A\rightarrow B$ is an integral extension of rings.
Let $\lp\in X$ be a point with $\lq=f(\lp)$. 

\begin{definition} We say $f$ is \emph{unramified at $\lp$} if the following two conditions holds:
The canonical map
\begin{align} 
&\label{c1} f^{\#}:\O_{Y,\lq}\rightarrow f_*\O_{X,\lp}\text{ satisfies } \lm_{\lq}\O_{X,\lp}=\lm_{\lp}.\\
&\label{c2}\text{ The field extension } \kappa(\lq)\subseteq
\kappa(\lp)\text{ is a finite separable extension.}
\end{align}
We say the morphism $f$ is \emph{etale at $\lp$} if it unramified at $\lp$
and the ring homomorphism
\begin{align} 
&\label{c3}f^{\#}:\O_{Y,\lq}\rightarrow f_*\O_{X,\lp}
\end{align}
is flat. The morphism $f:X\rightarrow Y$ is an
\emph{etale morphism} if it is etale at $\lp$ for all $\lp\in X$.
The morphism $f$ is \emph{ramified} if for all $\lq\in Y$  there is a  $\lp \in f^{-1}(\lq)$
where $f$ is not etale at $\lp$.
\end{definition}

The following result is well known.
\begin{lemma} \label{etalecriteria1} The morphism $f$ is etale if and only if it is flat and
  for every $\lq\in Y$ the fiber $f^{-1}(\lq)$ is the disjoint union
  of reduced points $\lp\in f^{-1}(\lq)$ with $\kappa(\lq)\subseteq
  \kappa(\lp)$ a finite separable extension of fields.
\end{lemma}
\begin{proof} See \cite{groth} or \cite{milne}
\end{proof}

For a finite etale morphism $f:\Spec(B)\rightarrow \Spec(A)$ it is well
known the number of points $\lp$ in $f^{-1}(\lq)$ is constant. We let
$d=\# f^{-1}(\lq)$ be the \emph{degree} of $f$.

Recall the following general result:
Let $A_1,..,A_l$ be commutative rings with unit and let $S_i\subseteq
A_i$ be multiplicatively closed subsets for $i=1,..,l$. Let
$A=\oplus_{i=1}^l A_i$ be the direct sum of the commutative rings $A_i$ and
let $S=\oplus_{i=1}^l S_i$. It follows $S\subseteq A$ is a multiplicatively
closed subset.

\begin{lemma} There is an isomorphism\label{localization}
\[ (S_1\oplus \cdots \oplus S_l)^{-1}(A_1\oplus \cdots \oplus A_l) 
\cong  S_1^{-1}A_1\oplus \cdots \oplus S_l^{-1}A_l \]
of commutative rings.
\end{lemma}
\begin{proof} We prove this by induction on $l$. Assume $l=2$. We want
  to prove the ismorphism
\[ (S\oplus T)^{-1}(A\oplus B) \cong S^{-1}A\oplus
T^{-1}B.\]
Define the following morphism
\[ g: A \oplus B \rightarrow S^{-1}A\oplus T^{-1}B \]
by
\[ g(a,b)=(a/1,b/1).\]
It follows for all $(s,t)\in S\oplus T$ the element
$g(s,t)$ is invertible. Moreover if $g(a,b)=(a/1,b/1)=0$ it
follows $a/1=0=b/1$ hence there is an element $(s,t)\in S\oplus T$
with $sa=tb=0$. It follows $(s,t)(a,b)=(sa,tb)=0$. Finally any element
$(a/s,b/t)\in S^{-1}A\oplus T^{-1}B$ may be written as
\[ g(a,b)g(s,t)^{-1} .\]
It follows there is a canonical isomorphism
\[ S^{-1}A\oplus T^{-1}B \cong (S\oplus T)^{-1}(A\oplus
B) \]
of rings and the claim is proved. The Lemma now follows by induction.
\end{proof}

Let $K$ be an arbitrary field and let $P(t)=t^d+a_1t^{d-1}+\cdots
+a_{d-1}t+a_d\in K[t]$ be a polynomial with coefficients in
$K$. Recall the following notion: The polynomial $P(t)$ is
\emph{separable} if its roots in the algebraic closure  $\overline{K}$
of $K$ are all distinct. The polynomial $P(t)$ is \emph{inseparable}
if it has multiple roots.
Recall the following well known result:
Let $X=\Spec(K[t]/P(t)$ and  $Y=\Spec(K)$. Let $f:X\rightarrow Y$ be the
structure morphism.

\begin{proposition} \label{separable} The morphism $f$ is etale if and only if $P(t)\in K[t]$
  is a separable polynomial. The morphism $f$ is ramified if and only if $P(t)$ is inseparable.
\end{proposition}
\begin{proof}
Assume $P(t)$ is separable. It follows $P(t)$ has $d$ distinct roots
$\alpha_1,..,\alpha_d\in \overline{K}$. Let $P(t)=Q_1^{q_1}\cdots
Q_m^{q_m}$ be a decomposition
of $P(t)$ in $K[t]$ where $Q_i$ are irreducible polynomials in
$K[t]$. Since $P(t)$ is a separable polynomial it follows $q_i=1$ for
$i=1,..,m$. Let
$X_i=\Spec(K[t]/Q_i(t))$. It follows by the Chinese Remainder Theorem
that
\[ K[t]/P(t)\cong K[t]/Q_1(t) \oplus \cdots \oplus
K[t]/Q_m(t) .\]
We get an isomorphism
\[ X\cong X_1\cup \cdots \cup X_m \]
of schemes. We get a map
\[ f:X_1\cup \cdots \cup X_m \rightarrow Y \]
induced by the natural map $K\rightarrow K[t]/P(t)$. 
Let $L=K[t]/P(t)$ and $L_i=K[t]/Q_i(t)$. It follows $L_i$ is a
separable field extension of $K$. By the Kunneth formula the following
holds: There is an isomorphism
\[ \H^0(X,\O_X)\cong \H^0(X_1,\O_{X_1})\oplus \cdots \oplus
\H^0(X_m,\O_{X_m}) \]
of rings. One also sees $L_i=\H^0(X_i,\O_{X_i})$.
Let $\m_i$ be the following ideal:
\[ \m_i=L_1\oplus \cdots \oplus L_{i-1}\oplus \{0\} \oplus
L_{i+1}\oplus \cdots \oplus L_m.\]
It follows 
\[ \H^0(X,\O_X)/\m_i\cong \H^0(X_i,\O_{X_i})=L_i \]
hence $\m_i\subseteq \H^0(X,\O_X)$ is a maximal ideal. The ideals
$\m_1,..,\m_m$ are all maximal ideals in $\H^0(X,\O_X)$.
Consider $\m_i\in X$ and look at the map
\[ f^{\#}:\O_{Y,f(\m_i)}\rightarrow \O_{X,\m_i}.\]
It is given by the natural map $K\rightarrow L_{\m_i}$. Let
$S_i\subseteq \H^0(X,\O_X)$ be the multiplicatively closed subset
defined by $S_i=\H^0(X,\O_X)-\m_i$. It follows
\[ S_i=L_1\oplus \cdots \oplus L_i^*\oplus \cdots \oplus L_m \] hence
we get by lemma \ref{localization} an isomorphism
\[ L_{\m_i}\cong S_i^{-1}\H^0(X,\O_X)\cong \]
\[ L_1^{-1}L_1\oplus \cdots \oplus (L_i^*)^{-1}L_i \oplus \cdots \oplus
L_m^{-1}L_m \cong L_i.\]
Hence the map
\[ f^{\#}:\O_{Y,f(\m_i)}\rightarrow \O_{X,\m_i} \]
is the map
\[ K\rightarrow L_i .\]
The field extension
\[ \kappa(f(\m_i))=K\subseteq L_i=\kappa(\m_i) \]
is separable.
It follows the map $f$ is flat and
unramified at $\m_i$ hence $f$ is etale.
Conversely assume $f:X\rightarrow Y$ is etale and
$P(t)$ inseparable. It follows $P(t)=Q_1^{q_1}\cdots Q_m^{q_m}$ where
all polynomials $Q_i$ are irreducibel and one of the following holds:
\begin{align}
&\text{There is an $i$ with $q_i>1$. }\\
&\text{All $q_i=1$ and there is an $i$ where $K\subseteq K[t]/Q_i$ is
  inseparable.}
\end{align}
Let $\lm_i$ be the maximal ideal corresponding to $i$. It follows $f$
is not etale at $i$ which is a contradiction, and the first part of
the Proposition is proved.
The second part of the claim  is obvious and the Proposition is proved.
\end{proof}

Let 
\[ F(t)=a_mt^m+a_{m-1}t^{m-1}+\cdots +a_1t+a_0 \]
and
\[ G(t)=b_nt^n+b_{n-1}t^{n-1}+\cdots +b_1t+b_0 \]
be polynomials in $A[t]$ where $A$ is an arbitrary commutative ring
with unit and $a_m,b_n\neq 0$. Make the following definition
\begin{definition} \label{resultant} Let $Res_{m,n}(F(t),G(t))$ be the following
  determinant
\[
\begin{vmatrix} a_m & a_{m-1}  & \cdots & a_1 & a_0 & 0 & 0 &
\cdots & 0 & 0  \\
              0 &  a_m & a_{m-1} & \cdots & a_1 & a_0 & 0 &   \cdots & 0
              & 0  \\
              0 & \vdots & \cdots &  \cdots & \cdots & \cdots & \cdots
              & \vdots & a_1 & a_0 \\
              b_n & b_{n-1}  & \cdots & b_1 & b_0 & 0 & 0 &
\cdots & 0 & 0  \\
              0 &  b_n & b_{n-1} & \cdots & b_1 & b_0 & 0 &   \cdots
              & 0
              & 0  \\
              0 & \vdots & \cdots &  \cdots & \cdots & \cdots & \cdots
              & \vdots & b_1 & b_0 \\
\end{vmatrix}.
\]
We say the element $Res_{m,n}(F(t),G(t))$ is the \emph{resultant}
of the polynomials $F(t)$ and $G(t)$.
\end{definition}
The matrix contains $n$ rows with $a_i$'s and $m$ rows with
$b_j$'s. It is a square matrix of rank $m+n$. We often write
$Res(F,G)$ instead of $Res_{m,n}(F(t),G(t))$. From the definition it
is immediate $Res(F,G)$ is an element of the ring $A$.

The following result is well known:
\begin{proposition} \label{sylvester}Let $A=K$ be an arbitrary field, and let
  $F(t),G(t)$ be polynomials in $K[t]$ of degree $m$ and $n$ with
  $a_m,b_n\neq 0$. Let $\alpha_1,..,\alpha_m$ and $\beta_1,..,\beta_n$
  be the roots of $F$ and $G$ in an algebraic closure $\overline{K}$
  of $K$. Assume $F'(t)$ is a polynomial of degree $m'\leq m$.
The following holds:
\begin{align}
&\label{s1}Res(F,G)=a_m^nb_n^m\prod_{i,j}(\alpha_i-\beta_j) \\
&\label{s2}Res(F,G)=0 \iff\text{ $F$ and $G$ have a common root in $\overline{K}$}\\
&\label{s3}Res(F,G)=(-1)^{mn}Res(G,F) \\
&\label{s4}Res_{m,n}(FF',G)=Res_{m,n}(F,G)Res_{m',n}(F',G)
\end{align}
\end{proposition}
\begin{proof} For a proof of the facts \ref{s1}-\ref{s4} see \cite{gelfand}, Section 12.
\end{proof}

Let $A$ be an arbitrary commutative ring with unit and let 
$P(t)=t^d+a_{d-1}t^{d-1}+\cdots +a_{1}t+a_0 \in A[t]$ be any degree d
monic polynomial. The formal derivative $P'(t)$ is again a polynomial
in $A[t]$.

\begin{definition} \label{discrdef} We let $Discr(P(t))=Res_{d,d-1}(P(t),P'(t))\in A$ be the 
\emph{discriminant} of the polynomial $P(t)$. We say $P(t)$ is
\emph{separable} if $Discr(P(t))$ is a unit in $A$. We say $P(t)$ is
\emph{inseparable} if $Discr(P(t))$ is nilpotent.
\end{definition}

\begin{proposition} \label{discriminant} Assume $A=K$ is a field. It follows
  $Discr(P(t))=0$ if and only if $P(t)$ has a root $\alpha\in
  \overline{K}$ of multiplicity greater than $2$.
\end{proposition}
\begin{proof} Assume $Discr(P(t))=0$ it follows $P(t)$ and $P'(t)$
  have a common root $\alpha \in \overline{K}$. It follows
  $P(t)=(t-\alpha)^2Q(t)$ for some polynomial $Q(t)\in \overline{K}[t]$,
  hence $\alpha$ has mutiplicity greater than $2$. Conversely, if
  $P(t)=(t-\alpha)^2Q(t)$ with $Q(t)\overline{K}[t]$ it follows $P(t)$
  and $P'(t)$ have a common root, hence $Discr(P(t))=0$ and the claim
  of the Proposition follows.
\end{proof}

\begin{corollary} \label{discr_criteria}The following holds:
  $Discr(P(t))\neq 0$ if and only if all roots of $P(t)$ 
have multiplicity one.
\end{corollary}
\begin{proof} Since $Discr(P(t))\neq 0$ it follows from Proposition
  \ref{discriminant} that all roots of $P(t)$ in $\overline{K}$ are of
  multiplicity one. Conversely
  if all roots of $P(t)$ have multiplicity one it follows $Discr(P(t))\neq 0$. The
  Corollary is proved.
\end{proof}

Let $\psi:A\rightarrow B$ be a map of
commutative rings. We get an induced map
\[ \psi_t:A[t]\rightarrow B[t] \]
defined by
\[\psi_t(P(t))=\psi(b_n)t^n+\psi(b_{n-1})t^{n-1}+\cdots +\psi(b_1)t+\psi(b_0).\]
Let $P_\psi(t)=\psi_t(P(t))$.

We get the following result:
\begin{lemma} There is an equality
\[ \psi(Discr(P(t))=Discr(P_\psi(t)) \]
in $B$.
\end{lemma}
\begin{proof}  By definition 
\[ P(t)=t^d+a_{d-1}t^{d-1}+a_{d-2}t^{d-2}+\cdots +a_1t+a_0 \]
has coefficients $a_i$ in $A$.
and
\[ P'(t)=dt^{d-1}+(d-1)a_{d-1}t^{d-2}+\cdots +2a_2t+a_1.\]
The discriminant $Discr(P(t),P'(t))$ is by Definition \ref{discrdef} given
by the following determinant:
\[|M|=
\begin{vmatrix} 1 & a_{d-1}  & \cdots & a_1 & a_0 & 0 & 0 & 
\cdots & 0 & 0  \\
              0 &  1 & a_{d-1} & \cdots & a_1 & a_0 & 0 &   \cdots & 0
              & 0  \\
              0 & \vdots & \cdots &  \cdots & \cdots & \cdots & \cdots
              & \vdots & a_1 & a_0 \\
              d & (d-1)a_{d-1}  & \cdots & a_1 & 0 & 0 & 0 &
\cdots & 0 & 0  \\
              0 &  d & (d-1)a_{d-1} & \cdots & a_1 & 0 & 0 &   \cdots & 0
              & 0  \\
              0 & \vdots & \cdots &  \cdots & \cdots & \cdots & \cdots
              & \vdots & 2a_2 & a_1 \\
\end{vmatrix}.
\]
Assume $\psi:A\rightarrow B$ is a ring homomorphism. 
It follows
\[ \psi(Discr(P(t))= \psi(|M|)=|\psi(M)|=Discr(P_\psi(t))\]
and the Lemma is proved.
\end{proof}

Let $b=Discr(P(t))\in A$. Consider the natural morphism
\[ \pi:\Spec(A[t]/P(t))\rightarrow \Spec(A) \]
of affine schemes. We get a diagram of maps of affine schemes
\[
\diagram  \pi^{-1}(V(b)) \rto^i \dto^{\pi} & \Spec(A[t]/P(t)) \dto^\pi
& \pi^{-1}(D(b))\lto^l  \dto^\pi \\
V(b)\rto^j & \Spec(A) & D(b)\lto^k 
\enddiagram
\]
where $i,j,k$ and $l$ are the natural inclusions.

\begin{proposition} \label{etale} Let $S=\{1,b,b^2,..\}\subseteq A$
  and let $A_b=S^{-1}A$. Let $\phi:A[t]\rightarrow A_b[t]$ be the
  natural map. The following holds:
\begin{align}
&\label{et1}\text{The morphism }\pi:\pi^{-1}(D(b))\rightarrow D(b)\text{ is
  etale.}\\
&\label{et2}\text{The morphism }\pi:\pi^{-1}(V(b))\rightarrow V(b)\text{ is ramified.}\\
&\label{et3}\text{There is an isomorphism }\pi^{-1}(D(b))\cong \Spec(A_b[t]/P_\phi(t))
\end{align}
\end{proposition}
\begin{proof} We first prove \ref{et1}: We want to show 
\[ \pi:\pi^{-1}(D(b))\rightarrow D(b) \]
is an etale morphism.
Pick $\lp \in \Spec(A)$ with $\lp\in D(b)$. It follows $b\notin \lp$. 
Let $\psi:A\rightarrow \kappa(\lp)$ be the natural map where
$\kappa(\lp)$ is the residue field of $\lp$. The induced map on the
fiber $\pi^{-1}(\lp)\rightarrow \Spec(\kappa(\lp))$ is the natural map
\[ \Spec(\kappa(\lp)[t]/P_\psi(t))\rightarrow \Spec(\kappa(\lp)) .\]
Since $b=Discr(P(t))\notin \lp$ it follows  
\[ \psi(Discr(P(t)))=Discr(P_\psi(t))\neq 0 \]
in the residue field $\kappa(\lp)$. It follows from Corollary
\ref{discr_criteria} $P_\psi(t)$ is separable in $\kappa(\lp)[t]$
hence the map
\[ \pi:\Spec(\kappa(\lp)[t]/P_\psi(t))\rightarrow \Spec(\kappa(\lp))
\]
is by Proposition \ref{separable} an etale map. It follows from Lemma
\ref{etalecriteria1} the map
\[ \pi:\pi^{-1}(D(b)) \rightarrow D(b) \]
is etale and claim \ref{et1} is proved.
We prove claim \ref{et2}: Pick a point $\lp \in V(b)$ and consider the
morphism
\[ \psi:\pi^{-1}(\lp)=\Spec(\kappa(\lp)[t]/P_\psi(t))\rightarrow
\Spec(\kappa(\lp)).\]
Since $\psi(b)\in \lp$ it follows $Discr(P_\psi(t))=0$ in the
residue field $\kappa(\lp)$. By  Proposition \ref{discriminant} it
follows  $P_\psi(t)$ is inseparable over
$\kappa(\lp)$. It follows from Proposition \ref{separable} for each $\lp\in V(b)$ the morphism
\[ \pi^{-1}(\lp)\rightarrow \Spec(\kappa(\lp)) \]
is ramified and claim \ref{et2} follows. 
We prove claim \ref{et3}: There is an isomorphism of rings 
\[ S^{-1}(A[t]/P(t))\cong S^{-1}(A[t])/P_\phi(t)\cong A_b[t]/P_\phi(t)
.\]
Hence 
\[ \pi^{-1}(D(b))\cong \Spec(S^{-1}(A[t]/P(t)))\cong
\Spec(A_b[t]/P_\phi(t)) .\]
Moreover, the natural map
\[\Spec(A_b[t]/P_\phi(t))\rightarrow \Spec(A_b) \]
is the map
\[ \pi^{-1}(D(b))\rightarrow D(b) \]
and the Proposition is proved.
\end{proof}

\begin{example} An inseparable polynomial.
 \end{example}
Let $a_1,..,a_d\in A$ be elements with $a_i-a_j$ nilpotent for some
$i\neq j$. It follows the polynomial
\[P(t)=(t-a_1)(t-a_2)\cdots (t-a_d)\in A[t] \]
is an inseparable polynomial: Let $\lp\subseteq A$ be a prime
ideal. It follows $a_i-a_j \in \lp$ hence $a_i=a_j$ in the residue
field $\kappa(\lp)$. Let $\phi:A\rightarrow \kappa(\lp)$ be the
canonical map. It follows $P_\phi(t)$ has multiple roots hence
$Discr(P_\phi(t))=\phi(Discr(P(t))=0$ in $\kappa(\lp)$. It follows
$Discr(P(t))$ is nilpotent since $Discr(P(t))\in \lp$ for all primes
$\lp$. It follows $P(t)$ is an inseparable polynomial.

\begin{example} A separable polynomial. \end{example}
Let $a_1,..,a_d\in A$ be elements with $a_i-a_j$ not nilpotent for all
$i\neq j$. It follows the polynomial
\[ P(t)=(t-a_1)(t-a_2)\cdots (t-a_d) \in A[t] \]
is a separable polynomial: Let $\lp \subseteq A$ be a prime ideal.
It follows $a_i\neq a_j $ in $\kappa(\lp)$ for all $i\neq j$. Let 
$\phi:A\rightarrow \kappa(\lp)$ be the canonical map. It follows
$P_\phi(t)$ is a separable polynomial for all $\lp$. It follows
$Discr(P_\phi(t))=\phi(Discr(P(t)))\neq 0$ in $\kappa(\lp)$ for all
primes $\lp$. Hence $Discr(P(t))$ is a unit in $A$ and $P(t)$ is a
separable polynomial.

\begin{corollary} \label{maximal}The open set $U=D(b)\subseteq \Spec(A)$
is the maximal open subset $U\subseteq \Spec(A)$ where $\pi:\pi^{-1}(U)\rightarrow U$
is etale.
\end{corollary}
\begin{proof} This follows from Proposition \ref{etale},
  Claim \ref{et1} and \ref{et2}.
\end{proof}

\begin{corollary} \label{criteria1} The morphism $\pi:\Spec(A[t]/P(t))\rightarrow
  \Spec(A)$ is etale if and only if $Discr(P(t))$ is a unit in $A$.
\end{corollary}
\begin{proof} By Proposition \ref{etale}, \ref{et2} it follows the morphism $\pi$ is etale if and only if
  $V(Discr(P(t)))=\emptyset$. This is if and only if $Discr(P(t))$ is a
  unit, and the Corollary follows.
\end{proof}

\begin{corollary} The morphism $\pi:\Spec(A[t]/P(t))\rightarrow
  \Spec(A)$ is ramified if and only if $Discr(P(t))$ is nilpotent in $A$.
\end{corollary}
\begin{proof} Let $\lp\in \Spec(A)$ and let $\phi:A\rightarrow
  \kappa(\lp)$ be the canonical map. It follows $\pi$ is ramified 
if and only if $P_\phi(t)\in \kappa(\lp)[t]$ is inseparable for all
primes $\lp$. This is
if and only if $Discr(P_\phi(t))=\phi(Disr(P(t))$ is zero in
the residue field $\kappa(\lp)$ for all primes $\lp$. This is if and only if $Disc(P(t))\in \lp$ for all
primes $\lp \subseteq A$. This is if and only if $Discr(P(t))$ is
nilpotent in $A$, and the Corollary follows.
\end{proof}

It follows the morphism $\pi$ is etale if and only if $P(t)$ is
separable. The morphism $\pi$ is ramified if and only if $P(t)$ is inseparable.
Hence the discriminant $Discr(P(t))$ measures when the morphism
\[ \pi:\Spec(A[t]/P(t))\rightarrow \Spec(A) \]
is etale.

\begin{definition} The scheme $D(\pi)=V(Discr(P(t))\subseteq \Spec(A)$
  is the \emph{discriminant of the morphism $\pi$} where
  $\pi:\Spec(A[t]/P(t))\rightarrow \Spec(A)$ is the canonical morphism.
\end{definition}

The underlying set of points of $D(\pi)$ equals the set theoretical
discriminant $Discr^1(\pi)$ from Definition \ref{settheoretical}.
By Proposition \ref{etale} $D(\pi)$ is the largest closed subscheme of
$\Spec(A)$ with the property that the map
 $\pi:\pi^{-1}(D(\pi))\rightarrow D(\pi)$  is ramified.

\begin{theorem} \label{standard}Let $f:X\rightarrow Y$ be an etale
  morphism of degree
  $d$ and let $\lp\in X$ with $\lq=f(\lp)$. There is an open affine
  neighborhood $U=\Spec(A)\subseteq Y$ with $\lq \in U$ and
  $f^{-1}(U)=\Spec(B)$ where $B=A[t]/P(t)$ where $P(t)$ is a monic
  degree $d$ polynomial with $Discr(P(t))$a unit in $A$.
\end{theorem}
\begin{proof} By \cite{groth} it follows there are affine open sets
$\lp\in U=\Spec(B)\subseteq X$ and $V=\Spec(A)\subseteq Y$ with
$f(U)\subseteq V$ and a commutative diagram
\[
\diagram U \dto^f \rto_j & \Spec(A[t]_{P'(t)}/P_\phi(t)) \dto^f \\
          V \rto^\cong & \Spec(A)
\enddiagram
\]
where $j$ is an open immersion, 
\[ \phi:A[t]\rightarrow A[t]_{P'(t)} \]
is the canonical map and $P(t)\in A[t]$ is a monic degree $d$
polynomial. Since $U\subseteq
\Spec(A[t]_{P'(t)}/P_\phi(t))=D(P'(t))$ it follows by Proposition
\ref{etale} $f(U)\subseteq D(b)$
where $b=Disrc(P(t))\in A$.
Since the extension $A\subseteq A[t]/P(t)$ is faitfully flat it
follows
\[ \Spec(A[t]/P(t))\rightarrow \Spec(A) \]
is open.
It follows the morphism $f:U\rightarrow V$ is
open hence $f(U)\subseteq V$ is an open set. We may therefore choose a
basic
open set $\Spec(A_a)\subseteq f(U)$ with $\lq=f(\lp)\in
\Spec(A_a)$. It follows the map
\[ f:f^{-1}(\Spec(A_a))\rightarrow \Spec(A_a) \]
equals the map
\[ f: \Spec(A_a[t]/P_\phi(t))\rightarrow \Spec(A_a) \]
where
\[ \phi:A[t]\rightarrow A_a[t] \]
is the canonical map. Since $\Spec(A_a)\subseteq D(b)$ it follows by
Proposition \ref{etale} and Corollary \ref{criteria1} 
$Discr(P_\phi(t))$ is a unit, and the claim of the Theorem follows.
\end{proof}

\begin{definition} \label{standardetale} A map
\[ \Spec(A[t]/P(t))\rightarrow \Spec(A) \]
where $P(t)\in A[t]$ is a monic polynomial with $Discr(P(t))$ a unit is
called a \emph{standard etale morphism}.
\end{definition}

By Theorem \ref{standard} it follows any etale morphism is locally a
given by a standard etale morphism.

Note: The definition of standard etale morphism given in Definition
\ref{standardetale} differs from the one give in the litterature (see
\cite{groth} and \cite{milne}). 
In the litterature a standard etale morphism is a morphism on the form
\[ \pi:\Spec(A[t]_{P'(t)}/P_\phi(t))\rightarrow \Spec(A) \]
where 
\[ \phi:A[t]\rightarrow A[t]_{P'(t)} \]
is the canonical map.

Let $P(t)\in A[t]$ be any monic degree $d$ plynomial where $A$ is an
arbitrary commutative ring. Let $b=Discr(P(t))\in A$ be the
discriminant. By Proposition \ref{etale} we get a diagram of maps of
schemes
\[
\diagram \pi^{-1}(V(b)) \rto^i \dto^\pi & \Spec(A[t]/P(t)) \dto^\pi &
\Spec(A_b[t]/P_\phi(t)) \lto^j \dto^\pi \\
           V(b) \rto^k & \Spec(A) & \Spec(A_b) \lto^l
\enddiagram
\]
where $\phi:A\rightarrow A_b$ is the canonical morphism. It follows by
Corollary \ref{criteria1} the morphism
\[ \Spec(A_b[t]/P_\phi(t)) \rightarrow \Spec(A_b) \]
is standard etale. Hence any morphism
\[ \Spec(A[t]/P(t)) \rightarrow \Spec(A) \]
decompose into a standard etale morphism and a ramified morphism. The
open set $D(b)=\Spec(A_b)\subseteq \Spec(A)$ is dense, hence the
morphism
\[ \Spec(A[t]/P(t))\rightarrow \Spec(A) \] 
is generically etale.

Let $P(t)=a_dt^d+a_{d-1}t^{d-1}+\cdots +a_1t+a_0\in A[t]$ be any
polynomial with coefficients in $A$. Assume $a_d\neq 0$.
Assume we are given  integers $d>d_1>d_2> \cdots >d_k \geq 0$.
Let $A_i=A/(a_d,a_{d_1},..,a_{d_i})$ where $a_{d_i}$ is the $d_i$'th
coefficient of $P(t)$.
and $A^i=(A_{i-1})_{a_{d_i}}$. Let $U_i=\Spec(A^i)$ and $V_i=\Spec(A_i)$.

\begin{theorem} \label{main1} There exists unique integers $d>d_1>d_2> \cdots >d_k
  \geq 0$ 
satisfying the following: We may write $\Spec(A)=V_0\cup U_0$
and for all $i=1,..,k$  $V_i=V_{i+1}\cup U_{i+1}$.
Moreover for all $i$ there is an element $b_i\in A^i$ giving a
disjoint union
\[ \Spec(A^i)=V(b_i)\cup D(b_i) \]
with the following properties: The natural morphism
\[ \pi: \Spec(A^i[t]/\overline{P}(t))\rightarrow \Spec(A^i) \]
satisfy
\begin{align}
&\label{d1}\pi: \pi^{-1}(D(b_i))\rightarrow D(b_i) \text{ is standard etale of degree $d_i$.}\\
&\label{d2} \pi: \pi^{-1}(V(b_i))\rightarrow V(b_i) \text{ is ramified.}
\end{align}
The open set $D(b_i)\subseteq \Spec(A^i)$ is the largest open
subscheme with property \ref{d1}.
\end{theorem}
\begin{proof} Consider 
\[ P(t)=a_dt^d+\cdots +a_1t+a_0 \in A[t].\]
We may write $\Spec(A)=V(a_d)\cup D(a_d)$ as a disjoint union.
It follows $V(a_d)=\Spec(A_0)$ and $D(a_d)=\Spec(A^0)$. On $D(a_d)$ it
follows the leading coefficient $a_d$ of $P(t)$ is invertible. Let 
\[ \phi:A \rightarrow A_{a_d}=A^0 \]
be the canonical map of rings. It follows $P_\phi(t)\in A_{a_d}[t]$ has
a unit as leading coefficient. 
Let $b_0=Discr(P_\phi(t))\in A^0$ it follows $\Spec(A^0)=V(b_0)\cup
D(b_0)$.
We get a canonical morphism
\[ \pi:\pi^{-1}(D(b_0))=\Spec(A^0[t]/P_\phi(t))\rightarrow \Spec(A^0)=D(b_0) \]
and by Proposition \ref{etale}, \ref{et1} it follows the morphism
\[ \pi:\pi^{-1}(D(b_0))\rightarrow D(b_0) \]
is standard etale of degree $d$.
Again by Proposition \ref{etale}, \ref{et2} it follows the morphism
\[\pi:\pi^{-1}(V(b_0))\rightarrow V(b_0) \]
is ramified. By Corollary \ref{maximal} it follows $D(b_0)\subseteq
\Spec(A^0)$ is the maximal open subscheme with this property. The
Theorem now follows by Proposition \ref{etale} and an induction.
\end{proof}

\section{Discriminants of linear systems on the projective line}

In this section we study the discriminant $D^l(\O(d))$ of the line
bundle $\O(d)$ on $\p^1_K$ where $K$ is any field.

Let $E=K\{e_0,e_1\}$ be a $K$-vector space of dimension two where $K$
is any field and let $E^*=K\{x_0,x_1\}$ be its dual. Let
$X=\Spec(K)$ and $\p=\p(E^*)=\p^1_\Z\times_\Z \Spec(K)$ be the projective
line over $K$. Cover $\p$ by the two standard open subsets $D(x_i)$
and let $t=\frac{x_0}{x_1}$ and $s=\frac{x_0}{x_1}$. Let $\O(d)$ be
the invertible sheaf corresponding to the graded $K[x_0,x_1]$-module
$K[x_0,x_1](d)$, and let $W=\H^0(\p,\O(d))$ be its global
sections. It follows $W=\sym^d(E^*)$ hence $W$ is a free $K$-module on
the global sections $s_i=x_0^{d-i}x_1$. Let $y_i=s_i^*$. It follows
$\p(W^*)$ is the proj of the graded ring $K[y_0,..,y_d]$ and there is
a canonical structure morphism
\[ \pi:\p(W^*)\rightarrow \Spec(K) .\]
The tautological sequence on $\p(W^*)$ is the morphism
\begin{equation}\label{tautological}
 0\rightarrow \O(-1)\rightarrow \pi^*W .
\end{equation}
It is the sheafification of the following morphism of graded
$K[y_0,..,y_d]$-modules:
\[ \alpha:K[y_0,..,y_d](-1)\rightarrow K[y_0,..,y_d]\otimes_K K\{s_0,..,s_d\}
\]
\[  \alpha(1)=\sum_{j=0}^d y_j\otimes s_j .\]
Fix an integer $0\leq i \leq d$ and let $U_i=D(y_i)$. Restrict \ref{tautological} to $U_i$ to get the
following map:
\[ \alpha|_{U_i}:\O(-1)|_{U_i}\rightarrow \pi^*W|_{U_i} \]
given by
\[ \alpha(1/y_i)=\frac{\sum_{j=0}^d y_j\otimes
  s_j}{y_i}=\sum_{j=0}^d\frac{y_j}{y_i}\otimes s_j.\]
Let $u_j=\frac{y_j}{y_i}$ for $j=0,..,d$. It follows $u_i=\frac{y_i}{y_i}=1$.
Consider the diagram
\[
\diagram \p(W^*)\times_A \p \rto^p \dto^q & \p \dto^u \\
         \p(W^*) \rto^\pi & \Spec(K) 
\enddiagram
\]
The $l$'th \emph{Taylor map} for $\O(d)$ is a map of locally free
$\O_\p$-modules
\begin{equation}\label{taylor}
 T^l:u^*W\rightarrow \Pr^l(\O(d))
\end{equation}
on $\p=\p(E^*)$. Pull the Taylor map \ref{taylor} and the tautological
sequence back to $Y=\p(W^*)\times_K \p$ to get the sequence
\[ \O(-1)_Y \rightarrow^{\alpha_Y} W_Y \rightarrow^{T^l_Y} \Pr^l(\O(d))_Y \]
and define $\phi^l(\O(d))=T^l_Y\circ \alpha_Y$. 

\begin{definition} The scheme $I^l(\O(d))=Z(\phi^l(\O(d))$ is the
  \emph{l'th incidence scheme of } $\O(d)$. The scheme $D^l(\O(d))=q(I^l(\O(d)))$
is the \emph{$l$'th discriminant of} $\O(d)$.
\end{definition}

For $l=0$ the Taylor morphism
is the \emph{evaluation map}
\[ T^0:\O_\p\otimes W  \rightarrow \O(d) \]
defined locally by
\[ T^0(U):\O_\p(U)\otimes W \rightarrow \O(d)(U) \]
\[ T^0(U)(a\otimes s)=as|_U.\]

Let $U_{ij}=D(y_i)\times D(x_j)\subseteq \p(W^*)\times_K \p$ for fixed $i,j$.
Let $j=0$ and $t=x_1/x_0$.
Restrict the map $\phi^0(\O(d))$ to $U_{i0}$ to get the following map
\[ K[u_0,..,u_d][t]\frac{1}{y_i}\rightarrow K[u_0,..,u_d][t]\otimes_K
K\{s_0,..,s_d\} \rightarrow K[u_0,..,u_d][t]x_0^d \]
defined by
\[ \phi^0(\O(d))(\frac{1}{y_i})=u_0+u_1t+u_2t^2+\cdots +t_i+\cdots
+u_dt^d =f(t) .\]
It follows
\[ I^0(\O(d))|_{U_{i0}}=Z(f(t))\subseteq \Spec(K[u_0,..,u_d][t])=U_{i0}.\]

We get a map
\[ q|_{U_{i0}}:I^0(\O(d))|_{U_{i0}}\rightarrow D(y_i) .\]
Let $Z_{ij}=I^0(\O(d))|_{U_{ij}}$.

We get an induced map
\[ q^{\#}:\O_{D(y_i)}\rightarrow q_*\O_{Z_{i0}} \]
defined by the natural map
\[ q^{\#}:K[u_0,..,u_d]\rightarrow K[u_0,..,u_d][t]/(f(t)) .\]
It follows $Ker(q^{\#})=(0)$ hence $\overline{q(Z_{i0})}=D(y_i)$ for
$i=0,..,d$. It follows $D^0(\O(d))|_{D(y_i)}$ is schematically dense
in $D(y_i)$.

Consider the map $\phi^0(\O(d))$ on the open set $U_{i1}=D(y_i)\times
D(x_1)$.
Let $v=x_0/x_1$. We get the following map
\[ K[u_0,..,u_d][v]\frac{1}{y_i}\rightarrow K[u_0,..,u_d][v]\otimes_K
K\{s_0,..,s_d\} \rightarrow K[u_0,..,u_d][v]x_1^d \]
defined by
\[ \phi^0(\O(d))(\frac{1}{y_i})=u_0v^d+u_1v^{d-1}+u_2v^{d-2}+\cdots +v^{d-i}+\cdots
+u_d =g(v) .\]
It follows
\[ I^0(\O(d))|_{U_{i1}}=Z(g(v))\subseteq \Spec(K[u_0,..,u_d][v])=U_{i1}.\]
We get a map
\[ q|_{U_{i1}}:I^0(\O(d))|_{U_{i1}}\rightarrow D(y_i) .\]
Let $Z_{i1}=I^0(\O(d))|_{U_{i1}}$.
We get an induced map
\[ q^{\#}:\O_{D(y_i)}\rightarrow q_*\O_{Z_{i1}} \]
defined by the natural map
\[ q^{\#}:K[u_0,..,u_d]\rightarrow K[u_0,..,u_d][v]/(g(v)) .\]
It follows $Ker(q^{\#})=(0)$ hence $\overline{q(Z_{i1})}=D(y_i)$ for
$i=0,..,d$. It follows $D^0(\O(d))|_{D(y_i)}$ is schematically dense
in $D(y_i)$.

\begin{proposition} For all $d\geq 1$ it follows $D^0(\O(d))=\p(W^*)$.
\end{proposition}
\begin{proof} By the above argument it follows
  $\overline{q(I^0(\O(d)))}=\p(W^*)$ and the Proposition is proved.
\end{proof}

For all $1\leq l \leq d$ we get on $\p^1_K$
an exact sequence of locally free sheaves
\begin{equation}\label{exct}
 0\rightarrow \Q_{k,d} \rightarrow W\otimes \O_{\p^1_k}\rightarrow
\Pr^l(\O(d)) \rightarrow 0
\end{equation}
where $rk(\Pr^l(\O(d)))=l+1$, $rk(Q_{k,d})=d-k$ and $rk(W\otimes
\O_{\p^1_K})=d+1$.
Dualize the sequence \ref{exct} to get the sequence
\[ 0 \rightarrow \Pr^l(\O(d))^*\rightarrow W^*\otimes
\O_{\p^1_K}\rightarrow \Q_{k,d}^* \rightarrow 0.\]
Take relative projective space bundle to get a closed immersion
\[ \p(\Q_{k,d}^*)\subseteq \p(W^*)\times_K \p^1_K .\]
It follows from \cite{maa4} that $\p(\Q_{k,d}^*)=I^k(\O(d))$ hence
$dim(I^k(\O(d)))=dim(\p(\Q_{k,d}^*))=d-k-1+1=d-k$.

\begin{lemma} \label{irr1}Assume $p:X\rightarrow Y$ is any morphism of irreducible
  schemes with $U\subseteq X$ and $V\subseteq Y$ open
  subschemes with $p(U)\subseteq V$. Assume the induced morphism
  $\tilde{p}:U\rightarrow V$ is integral. It follows $dim(X)=dim(Y)$.
\end{lemma}
\begin{proof} Since $p:U\rightarrow V$ is integral it follows
  $dim(U)=dim(V)$. We get since $U,V$ are dense in $X$ and $Y$
\[ dim (X)=dim(U)=dim(V)=dim(Y) \]
and the Lemma follows.
\end{proof}

\begin{lemma} \label{irr2} Assume $p:X\rightarrow Y$ is a proper morphism of
  irreducible schemes satisfying the conditions of Lemma \ref{irr1}. 
Assume $Z\subseteq X$ is a closed irreducible subscheme. It follows $dim(Z)=dim(p(Z))$.
\end{lemma}
\begin{proof} Let $U=\Spec(B)\subseteq X$ and $V=\Spec(A)\subseteq Y$ be
  open dense subschemes such that $p$ induce an integral morphism
\[ \tilde{p}:U\rightarrow V .\] 
It follows 
$U\cap Z\subseteq Z$ and $V\cap p(Z)\subseteq p(Z)$ are open dense
subsets. 
We get an induced morphism
\[ q:U\cap Z\rightarrow V\cap p(Z) \]
which is integral, hence $dim(U\cap Z)=dim(V\cap p(Z))$. It follows
$dim(Z)=dim(p(Z))$ and the Lemma follows.
\end{proof}

Consider again the diagram
\[
\diagram Y=\p(W^*)\times_A \p \rto^p \dto^q & \p \dto^u \\
         \p(W^*) \rto^\pi & \Spec(A)
\enddiagram
\]
when $A=K$ is a field and look at the sequence of locally free sheaves
\[ \phi^l(\O(d)): \O(-1)_Y\rightarrow W_Y \rightarrow \Pr^l(\O(d))_Y .\]
Consider the Taylor map on $\p^1_K$. 
Let $U_i=D(x_i)$ and $t=x_1/x_0$, $s=x_0/x_1$.
On $U_0$ we get the following calculation:
\[ T^l|_{U_0}:\O_{U_{0}}\otimes W \rightarrow \Pr^l(\O(d))|_{U_0} \]
looks as follows:
\[ T^l|_{U_0}:K[t]\{ x_0^{d-i}x_1\} \rightarrow K[t]\{ dt^j\otimes
x_0^d\} \]
with 
\[ T^l(x_0^{d-i}x_1^i)=T^l(t^ix_0^d)=(t+dt)^i\otimes x_0^d .\]
Let $U_{ij}=D(y_i)\times D(x_j)\subseteq \p(W^*)\times \p^1$.
Let $u_j=y_j/y_i$.
The composed morphism
\[ \phi^l(\O(d))|_{U_{i0}}:K[u_0,..,u_d][t]\frac{1}{y_i}\rightarrow
K[u_0,..,u_d][t]\otimes_K\{x_0^{d-i}x_1\}
\rightarrow K[u_0,..,u_d][t]\{dt^j\otimes x_0^d\} \]
is the following map:
\[ \phi^l(\O(d))|_{U_{i0}}(1/y_i)=(f(t),f'(t)/1!,..,f^{(l)}(t)/l!) \]
where
\[ f(t)=u_0+u_1t+\cdots + t^i+\cdots +u_dt^d .\]
It follows we get an equality of ideal sheaves on $U_{i0}$
\[ \I_{I^l(\O(d))}|_{U_{i0}}=\{ f(t), f'(t),..,f^{(l)}(t) \} .\]
Consider $D(x_1)\subseteq \p^1$ and let $s=x_0/x_1$.
The Taylor map looks as follows
\[ T^l:K[s]\otimes \{ x_0^{d-i}x_1^i\}\rightarrow K[s]\{ ds^j\otimes
x_0^d\} \]
with 
\[ T^l(x_0^{d-i}x_1^i)=T^l(s^ix_0^d)=(s+ds)^{d-i}\otimes x_0^d.\]
The composed morphism
\[ \phi^l(\O(d))|_{U_{i1}}:K[u_0,..,u_d][s]\frac{1}{y_i}\rightarrow
K[u_0,..,u_d][s]\otimes_K\{x_0^{d-i}x_1\}
\rightarrow \]
\[ K[u_0,..,u_d][s]\{ds^j\otimes x_1^d\} \]
is the following map:
\[ \phi^l(\O(d))|_{U_{i1}}(1/y_i)=(g(s),g'(s)/1!,..,g^{(l)}(s)/l!) \]
where
\[ g(s)=u_0s^d+u_1s^{d-1}+\cdots + s^{d-i}+\cdots +u_d .\]
We get an equality of ideal sheaves
\[ \I_{I^l(\O(d))}|_{U_{i1}}=\{ g(s), g'(s),..,g^{(l)}(s) \} .\]
Let $V_{ij}=I^l(\O(d))|_{U_{ij}}$ and consider the morphism
\[ q_{ij}:V_{ij}\rightarrow q(V_{ij})=D^l(\O(d))\cap D(y_i).\]
Let $j=0$. We get an induced map at rings
\[ q_{i0}^{\#}:K[u_0,..,u_d]\rightarrow
K[u_0,..,u_d][t]/(f(t),..,f^{(l)}(t)) .\]
Let $P_j=Res(f^{(j)}(t),f^{(j+1)}(t))$. It follows there is an
equality
\[ ker(q_{i0}^{\#})=(P_0,..,P_{l-1}) .\]
We get an equality of ideal sheaves
\[ \I_{D^l(\O(d))}|_{D(y_i)}=\{P_0,..,P_{l-1}\} .\]

Consider the morphism
\[q_{i1}^{\#}:K[u_0,..,u_d]\rightarrow
K[u_0,..,u_d][s]/(g(s),..,g^{(l)}(s)/l!).\]
Let $Q_j=Res(g^{(j)}(s),g^{(j+1)}(s))$ it follows $Q_j=P_j$. It
follows
\[ ker(q_{i1}^{\#})=(Q_0,..,Q_{l-1})=(P_0,..,P_{l-1}) \]
hence $\I_{D^l(\O(d))}$ is locally generated by $l$ elements.

\begin{theorem} \label{main2} The $l$'th discriminant $D^l(\O(d))$ is an irreducible
local complete intersection of dimension $d-l$.
\end{theorem}
\begin{proof} Since $I^l(\O(d))=\p(\Q_{k,d}^*)$ is irreducible it
  follows $D^l(\O(d))$ is irreducible.
Let 
\[ q^l:I^l(\O(d)) \rightarrow D^l(\O(d)) \]
be the morphism induced by the projection morphism. Consider $l=0$ and
the open set $U_{d,0}=D(y_d)\times D(x_0)\subseteq \p(W^*)\times_K
\p^1_k$. Let $V_{d,0}=U_{d,0}\cap I^0(\O(d))$.
We get an induced morphism
\[ q^0:V_{d,0}\rightarrow q(V_{d,0}) \]
of schemes. Let $W_{d,0}=q(V_{d,0})$. We get an induced morphism of
sheaves
\[ q^{\#}:\O_{W_{d,0}}\rightarrow q_*\O_{V_{d,0}} \]
given by
\[ q^{\#}: B=K[u_0,..,u_d]\rightarrow K[u_0,..,u_d][t]/(f(t))=B[t]/(f(t)) \]
where
\[ f(t)=u_0+u_1t+\cdots +u_{d-1}t^{d-1}+t^d \]
It follows the ring extension $B\subseteq B[t]/(f(t))$ is an integral
extension hence the morphism
\[ q^0:V_{d,0}\rightarrow q(V_{d,0}) \]
is an integral morphism. By Lemma \ref{irr2} it follows
$dim(D^l(\O(d)))=dim(I^l(\O(d)))=d-l$.
Since the ideal sheaf of $D^l(\O(d))$ is locally generated by $l$
elemets the Theorem follows.
\end{proof}

\section{Discriminants of linear systems on flag varieties}

Let in the following $F$ be a fixed basefield of characeristic zero and
let $V$ be an $N+1$-dimensional vector space over $F$. Let $\p(V^*)$
be the projective space parametrizing lines in $V$. The space $\p(V^*)$
has the following property: There is
a well defined left action of $\SL(V)$ on $\p(V^*)$ and this give an
isomorphism $\SL(V)/P\cong \p(V^*)$ where $P\subseteq \SL(V)$ is the
parabolic subgroup fixing a line $L$ in $V$. The quotient $\SL(V)/P$
is a \emph{geometric quotient} in the sense of \cite{mumford}
and there is an equivalence of categories
between the category of linear finite dimensional representations of
$P$ and the category of locally free $\O_{\p(V^*)}$-modules with an
$\SL(V)$-linearization. 
We want to study the discriminant $D^k(\O(d))$ and the complex from
\ref{total} in this situation using $\SL(E)$-modules.

Let $\p=\p(V^*)$ and let $\I\subseteq \p\times \p$ be the ideal
of the diagonal. let $p,q:\p\times \p\rightarrow \p$ be the projection
morphisms.

\begin{definition} Let $k\geq 1$ and $d$ be integers. We define the
  \emph{$k$'th order jet bundle} of $\O(d)$ as follows:
\[ \Pr^k(\O(d))=p_*(\O_{\p\times \p}/\I^{k+1}\otimes q^*\O(d)).\]
\end{definition}

The invertible sheaf $\O(d)$ has a unique $\SL(V)$-linearization and
by prolongation of this structure it follows $\Pr^k(\O(d))$ has a
canonical $\SL(V)$-linearization. Hence $\Pr^k(\O(d))$ corresponds to
a unique finite dimensional $P$-module. The exterior product
$\wedge^j\Pr^k(\O(d))$ has for all $j\geq 1$ a canonical $\SL(V)$-linearization.
It follows the higher
cohomology group $\H^i(\p, \wedge^j\Pr^k(\O(d))  )$ is a finite
dimensional $\SL(V)$-module for all $i\geq 0$.

\begin{theorem} \label{exterior}There is for any integers $1\leq k <d$ and $1\leq
  j\leq  rk(\Pr^k)$  an isomorphism
\[ \H^0(\p, \wedge^j \Pr^k(\O(d))) \cong \sym^{j(d-k)}(V^*)\otimes
\wedge^j\sym^k(V^*) \]
of $\SL(V)$-modules. There is an equality $\H^i(\p,
\wedge^j\Pr^k(\O(d)))=0$ if $i>0$.
There is an isomorphism of $\SL(V)$-modules
\[ \H^i(\p, \wedge^j\Pr^k(\O(d))^*)= \sym^{j(d-k)-n-1}(V)\otimes
\wedge^j\sym^k(V^*) \]
if $i=n$ and $j(d-k)-n-1\geq 0$.
If $i=0,..,n-1$ or $i=n$ and $j(d-k)-n-1<0$ it follows $\H^i(\p,
\wedge^j\Pr^k(\O(d))^*)=0$.
\end{theorem}
\begin{proof} Let $\pi:\p\rightarrow S=\Spec(F)$ be the structure
  morphism. The $\SL(V)$-module $\sym^k(V^*)$ is a locally free
  $\O_S$-module with an $\SL(V)$-linearization. Pull back preserves
  the $\SL(V)$-linearization hence $\pi^*\sym^k(V^*)$ is a locally
  free  $\O_\p$-module with an $\SL(V)$-linearization. We may consider
  the following locally free sheaf:
\[ \O(d-k)\otimes \pi^*\sym^k(V^*).\]
Its corresponding $P$-module is $\sym^k(V^*)\otimes \sym^{d-k}(L^*)$
hence by \cite{maa1}, Theorem 2.4 we get an isomorphism of locally
free sheaves with $\SL(V)$-linearization
\[ \O(d-k)\otimes \pi^*\sym^k(V^*)\cong \Pr^k(\O(d)) .\]
We get the following calculation:
\[ \wedge^j \Pr^k(\O(d))\cong \wedge^j(\O(d-k)\otimes
\pi^*\sym^k(V^*))\cong \]
\[ \O(j(d-k))\otimes \wedge^j\pi^*\sym^k(V^*) \cong \O(j(d-k))\otimes
\pi^*\wedge^j\sym^k(V^*)\]
of $\SL(V)$-bundles.
It follows there is for every  $j\geq 1$ an isomorphism
\[ \wedge^j \Pr^k(\O(d))\cong \O(j(d-k))\otimes\pi^* \wedge^j\sym^k(V^*) \]
of $\SL(V)$-bundles.
We get the following calculation:
\[ \H^i(\p, \wedge^j \Pr^k(\O(d)))=\R^i\pi_*(\O(j(d-k))\otimes
\pi^*\wedge^j\sym^k(V^*))\cong\]
\[\wedge^j\sym^k(V^*)\otimes \R^i\pi_*\O(j(d-k)) .\]
The first part of the Theorem now follows from the calculation of equivariant cohomology of
invertible sheaves on projective space (see \cite{jantzen}).

By Theorem 2.4 in \cite{maa1} We get an isomorphism
\[ \wedge^j\Pr^k(\O(d))^*\cong \O(j(k-d))\otimes\pi^*\wedge^j
\sym^k(V) \]
as $\SL(V)$-bundles. We get using equivariant higher direct images the
following calculation:
\[ \H^i(\p,
\wedge^j\Pr^k(\O(d))^*)=\R^i\pi_*(\O(j(k-d))\otimes\pi^*\wedge^j\sym^k(V))=\]
\[ \wedge^j\sym^k(V)\otimes
\R^i\pi_*(\O(j(k-d))=\wedge^j\sym^k(V)\otimes \H^i(\p, \O(j(k-d)).\]
The second part of the Theorem  now follows from the calculation of equivariant
cohomology of invertible sheaves on projective space (see \cite{jantzen}).
\end{proof}

In several papers (see \cite{maa1}, \cite{piene} and \cite{dirocco})
the structure of the jet bundle on the projective line and projective
space has been studied. In the paper \cite{maa1} the $P$-module
structure of the jet bundle on projective space was classified. The
novelty of the result in Theorem \ref{exterior} is the calculation of
the $\SL(V)$-module structure of the higher cohomology groups
\[ \H(\P(V^*), \wedge^j\Pr^k(\O(d)) ) .\]
This is as indicated in the proof above a 
consequence of the result given in \cite{maa1}, Theorem 2.4 and
equivariant projection formulas.

\begin{example} Digression: The Borel-Weil-Bott Theorem .\end{example}
For all invertible sheaves $\O(d)$ on $\p$ with $d\geq 1$ the $\SL(V)$-module $\H^0(\p, \O(d))=\sym^d(V^*)$
is irreducible. This is a particular case of the \emph{Borel-Weil-Bott
Theorem}. We see from Theorem \ref{exterior} that this is no longer
true if we consider higher rank $\SL(V)$-linearized locally free
sheaves on $\p$. The $\SL(V)$-module
\[ \H^0(\p,  \wedge^j\Pr^k(\O(d)) )=\sym^{j(d-k)}(V^*)\otimes \wedge^j
\sym^k(V^*) \]
is never irreducible. 

Let $H\subseteq G$ be a closed subgroup.
We refer to a $G$-module of the form $\H^i(G/H,
\E(\rho))$ where $i\geq 0$ and $\E(\rho)$ is a $G$-linearized locally free
$\O_{G/H}$-module as a \emph{geometric $G$-module}.

The following general result is true: Let $P\subseteq G$ be a parabolic
subgroup of a linear algebraic group of finite type over $F$. The
quotient $G/P$ is a smooth projective scheme of finite type over
$F$ and $G$-linearized locally free $\O_{G/P}$-modules are in one to
one correspondence with rational  $P$-modules.

\begin{proposition} \label{global} Let $W$ be any $G$-module. It
  follows $W$ is a geometric $G$-module.
\end{proposition}
\begin{proof} Let 
\[\pi:G/P\rightarrow \Spec(F) \]
be the structure morphism. Let $\rho:G\rightarrow \GL(W)$ and 
consider the $\O_{G/P}$-module  $\O_{G/P}\otimes
\pi^*W$. It follows 
\[ (\O_{G/P}\otimes \pi^*W)(e)\cong F\otimes_F W\cong W \]
hence $\O_{G/P}\otimes \pi^*W$ is the $G$-linearized locally free
sheaf $\E(\rho)$ corresponding to $\rho$.
Since $W$ is a $G$-module it is also a $P$-module and the
locally free sheaf $\E(\rho)$ is a $G$-linearized  $\O_{G/P}$-module of rank
$r=dim(W)$. It is trivial as abstract locally free $\O_{G/P}$-module.
It follows from the equivariant projection formula (see \cite{jantzen})
\[ \H^0(G/P, \E(\rho))\cong \R^0\pi_*(\O_{G/P}\otimes \pi^*W) \cong \]
\[ W\otimes \R^0\pi_*\O_{G/P}\cong W\otimes_F F\cong W \]
since $\H^0(G/P, \O_{G/P})$ is the trivial rank one $G$-module, and the Proposition follows.
\end{proof}

Hence any $G$-module may be realized as the higher cohomology group of a
$G$-linearized locally free $\O_{G/P}$-module.

\begin{example} Decomposition of geometric $\SL(V)$-modules.\end{example}

One may ask the following general question:
Assume $G$  is a semisimple linear algebraic group and $P\subseteq G$
a parabolic subgroup.  Let $\L\in \Pic^G(G/P)$. One seeks a
decomposition
\[ \H^i(G/P, \wedge^j\Pr^k(\L))\cong \oplus_\lambda V_\lambda \]
of the $i$'th cohomology group of $\wedge^j\Pr^k(\L)$ into irreducible
$G$-modules $V_\lambda$. Since $G$ is semi simple and $\H^i(G/P,
\wedge^j\Pr^k(\L))$ is a finite dimensional $G$-module, such a
decomposition always exist by the general theory of representations of
semi simple algebraic groups. From Theorem \ref{exterior} one gets such
a decomposition on $\SL(V)/P=\p(V^*)$ by applying well known
combinatorial formulas from the theory of Schur-Weyl modules. 
The $\SL(V)$-module
\[ \H^0(\p, \wedge^j \Pr^k(\O(d))) \cong \sym^{j(d-k)}(V^*)\otimes
\wedge^j\sym^k(V^*) \]
is formed by applying compositions of Schur-Weyl modules to the
standard representation $V$ and its dual $V^*$. Given two Schur-Weyl
modules $\schur_\lambda$ and $\schur_\mu$ there is in general a
decomposition
\begin{equation} \label{plethysm}
 \schur_\lambda ( \schur_\mu(V^*))=\oplus_\lambda V_\lambda
\end{equation}
of $\schur_\lambda (\schur_\mu(V^*))$ into irreducible
$\SL(V)$-modules. To calculate this decomposition - referred to as \emph{plethysm} -
is an unsolved problem in general. There are several formulas
which are special cases of \ref{plethysm}: The \emph{Clebsch-Gordan formula}
which describes the decomposition of the tensor product of two
irreducible modules into irreducibles, 
and the \emph{Cauchy formula} which describe the decomposition of the
symmetric product of the tensor product of two standard modules (see \cite{fulton}). 
We leave it to the reader as an exercise to calculate this decomposition in the case of higher
cohomology groups of exterior powers of jet bundles on projective space.

Let $W=\H^0(\p,\O(d))$.
There is by \cite{maa1} on $\p$ an exact sequence of locally free sheaves
\[ 0\rightarrow Q_{k,d}\rightarrow \pi^*W \rightarrow
\Pr^k(\O(d))\rightarrow 0 \]
where the rightmost map is $T^k$ - the \emph{Taylor map}. 
It follows $\Q_{k,d}$ is a locally free sheaf.
Dualize this sequence to get the short exact sequence
\[ o\rightarrow \Pr^k(\O(d))^* \rightarrow \pi^*W^*
\rightarrow Q_{k,d}^*\rightarrow 0 .\]
Take relative projective space bundle to get the following closed
immersion of schemes:
\[ \p(Q_{k,d}^*)\subseteq \p(\pi^*W^*)\cong \p(W^*)\times \p  .\]

\begin{lemma} There is an equality 
\[ \p(\Q_{k,d}^*)=I^1(T^k) \]
of schemes.
\end{lemma}
\begin{proof}  Since $\Q_{k,d}^*=Coker((T^k)^*)$ the claim follows
from Theorem \ref{equality}.
\end{proof}

There is a commutative diagram of maps of schemes
\[
\diagram \p(Q_{k,d}^*) \rto^i \dto^{\tilde{q}} & \p(W^*)\times \p
\rto^p \dto^q & \p \dto^{\pi} \\
D^1(T^k) \rto^j & \p(W^*) \rto^{\pi} & S \\
\enddiagram \]

where $i,j$ are the inclusion morphisms. 

\begin{definition} The scheme $D^k(\O(d))=D^1(T^k)$ is the $k$'th discriminant
  of the invertible sheaf $\O(d))$. 
\end{definition}

\begin{example} \label{irreduciblepone} Irreducibility of $D^k(\O(d))$.\end{example}
Note: Since $\p(Q_{k,d}^*)$ is a projective bundle on $\p(V^*)$ it is
irreducible.  It  follows $D^k(\O(d))$ is
irreducible for all $1\leq k < d$. 

The morphism 
\[ \tilde{q}:\p(Q_{k,d}^*)\rightarrow D^k(\O(d)) \]
is by definition a surjective map of schemes. 
Let $Y=\p(W^*)\times \p$. On $\p(W^*)$ there is an exact sequence of locally free
sheaves
\[ 0 \rightarrow \O(-1) \rightarrow \pi^*W.\]
Pull this map back to $Y$ to get the composed map
\[ \psi:\O(-1)_Y\rightarrow W_Y \rightarrow ^{T^k}
\Pr^k(\O(d))_Y \]

\begin{proposition} There is an equality $Z(\psi)=\p(Q_{k,d}^*)$ of schemes.
\end{proposition}
\begin{proof} The claim follows from Theorem \ref{equality}  since $\Q_{k,d}^*=Coker((T^k)^*)$.
\end{proof}
We get an exact sequence
\[ \O(-1)_Y\otimes \Pr^k(\O(d))_Y^*\rightarrow \O_Y\rightarrow
\O_{Z(\psi)}\rightarrow 0.\]

The ideal sheaf of $Z(\psi)$ 
is locally generated by a regular sequence hence $Z(\psi)$ is a  local
complete intersection. 
There is by general results (see \cite{altman}) a Koszul complex of sheaves on Y:
\[ 0 \rightarrow \wedge^N\O(-1)_Y\otimes \Pr^k(\O(d))^*_Y \rightarrow \cdots \]
\[ \cdots \rightarrow \wedge^2\O(-1)_Y\otimes \Pr^k(\O(d))^*_Y
\rightarrow \O(-1)_Y\otimes \Pr^k(\O(d))^*_Y \rightarrow \]
\[ \O_Y\rightarrow \O_{Z(\psi)} \rightarrow 0.\]
Hence we get a resolution of the sheaf $\O_{Z(\psi)}$ on $Y$.
There is an isomorphism
\[ \Pr^k(\O(d))^*\cong \O(k-d)\otimes \pi^*\sym^k(V) \]
and an isomorphism
\[ \wedge^j\O(-1)_Y\otimes\Pr^k(\O(d))_Y\cong q^*\O(-j)\otimes
\wedge^jp^*(\O(k-d)\otimes \pi^*\sym^k(V))\cong \]
\[ q^*\O(-j)\otimes p^*\O(j(k-d))\otimes \pi^*\wedge^j\sym^k(V)\cong
\O(-j, j(k-d))\otimes \pi^*\wedge^j\sym^k(V).\] 
The complex now becomes
\[ 0\rightarrow \O(-N,N(k-d))\otimes \pi^*\wedge^N\sym^k(V) \rightarrow \cdots \]
\[\cdots \rightarrow \O(-j,j(k-d))\otimes \pi^*\wedge^j\sym^k(V)
\rightarrow \cdots \]
\[ \cdots \rightarrow \O(-1,k-d)\otimes \pi^*\sym^k(V)\rightarrow \O_Y
\rightarrow \O_{Z(\psi)}\rightarrow 0.\]

There is a projection morphism
\[ q:\p(W^*)\times \p \rightarrow \p(W^*) \]
and we want to push down the complex above to get a double complex on
$\p(W^*)$. When we push the complex down to $\p(W^*)$ we get a
double complex with terms given as follows:

\begin{theorem} \label{directimage} There is an isomorphism of $\SL(V)$-linearized locally
  free sheaves
\[ \R^nq_*(\wedge^j\O(-1)_Y\otimes \Pr^k(\O(d))^*_Y)\cong
\O(-j)\otimes \pi^* \sym^{j(d-k)-n-1}(V)\otimes \wedge^j\sym^k(V^*) \]
if $j(d-k)-n-1\geq 0$. If $i=0,..,n-1$ or $i=n$ and
$j(d-k)-n-1<0 $ it follows 
$\R^iq_*(\wedge^j\O(-1)_Y\otimes \Pr^k(\O(d))^*_Y)=0$.
There is an isomorphism of $\SL(V)$-linearized locally free sheaves
\[ \R^0q_*(\wedge^j\O(-1)_Y\otimes \Pr^k(\O(d))_Y\cong \O(-j)\otimes
\pi^*\sym^{j(d-k)}(V^*)\otimes \sym^k(V^*) .\]
If $i>0$ there is an equality $\R^iq_*(\wedge^j\O(-1)_Y\otimes \Pr^k(\O(d))_Y)=0$.
\end{theorem}
\begin{proof} We get by the projection formula and higher direct images 
for locally free sheaves with $\SL(V)$-linearization 
(see \cite{jantzen}) the following calculation:
\[ \R^iq_*(\wedge^j\O(-1)_Y\otimes \Pr^k(\O(d))^*_Y)\cong \]
\[ \R^iq_*(\wedge^j q^*\O(-1)\otimes p^*\Pr^k(\O(d))^*) \cong \]
\[\R^iq_*(q^*\O(-j)\otimes \wedge^j p^*\Pr^k(\O(d))^*)\cong \]
\[\O(-j)\otimes \R^iq_*p^*(\wedge^j \Pr^k(\O(d))^*)\cong \]
\[ \O(-j)\otimes \H^i(\p, \wedge^j\Pr^k(\O(d))^*) .\]
The first part of the Theorem now follows from Theorem \ref{exterior}.
We get the following calculation:
\[ \R^iq_*(\wedge^j\O(-1)_Y\otimes \Pr^k(\O(d))_Y)\cong \O(-j)\otimes
\pi^*\H^i(\p, \wedge^j\Pr^k(\O(d)) ).\]
The Theorem now follows from Theorem \ref{exterior}.
\end{proof}

\begin{example} \label{complex} Resolutions of ideal sheaves of
  discriminants.\end{example}

Assume $d-k-n-1\geq 0$, $l(j)=j(d-k)-n-1$ and $r=rk(\Pr^k(\O(d))$. 
The complex \ref{total} gives a complex of $\O_{\p(V^*)}$-modules
\[ 0\rightarrow \O(-r)\otimes \sym^{l(r)}(V)\otimes \wedge^r
\sym^k(V^*)  \rightarrow \cdots\]
\[ \cdots \rightarrow \O(-j)\otimes \sym^{l(j)}(V)\otimes
\wedge^j\sym^k(V^*) \rightarrow \cdots \]
\[ \O(-1)\otimes \sym^{l(1)}(V)\otimes \sym^k(V^*)\rightarrow
\O_{\p(V^*)}\rightarrow \O_{D^k(\O(d))}\rightarrow 0.\]
The hope is this complex can be used to construct a resolution of $D^k(\O(d))$.

\begin{example} \label{flag} Discriminants of linear systems on flag varieties.
\end{example}

In the following we use the notation of \cite{maa3}. Let $G=\SL(E)$
where $E$ is a
vector space of finite dimension over a field $F$ of characteristic
zero.
Let
\[ E_\bullet:0\neq E_1\subseteq E_2\subseteq \cdots \subseteq
E_k\subseteq E_{k+1}=E \]
be a flag in $E$ of type $\underline{d}=\{d_1,..,d_k\}$. Let
$n_i=d_1+\cdots +d_i$. It follows $dim(E_i)=n_i$. Let $P\subseteq
G=\SL(E)$ be the subgroup fixing the flag $E_\bullet$. It follows $P$
is a parabolic subgroup and the quotient $G/P$ is the flag
variety of $E$ of type $\underline{d}$. Let $V_\lambda$ be an
irreducible
$\SL(E)$-module with highest weight
\[\lambda=\sum_{i=1}^k l_i(L_1+\cdots +L_{n_i})=\sum_{i=1}^k l_i\omega_{n_i}.\]
Let $\L(\underline{l})\in \Pic^{G}(G/P)$ be the line bundle
corresponding to $\underline{l}=\{l_1,..,l_k\}\in\Z^k=\Pic^{G}(G/P)$.
It follows there is by Theorem 2.2 in \cite{maa3} an isomorphism
\[ V_\lambda \cong \H^0(G/P, \L(\underline{l}))^* \]
of $G$-modules.
Let
\[ T^k:\H^0(G/P, \L(\underline{l}))\otimes \O_{G/P} \rightarrow
\Pr^k(\L(\underline{l})) \]
be the Taylor map of order $k$.

\begin{theorem} \label{irr} The discriminant
  $D^k(\L(\underline{l}))=D^1(T^k)$ is irreducible
for all $1\leq k \leq min\{l_i+1\}$.
\end{theorem}
\begin{proof}
From \cite{maa3}, Theorem 3.7 it follows the Taylor map
\[ T^k:\H^0(G/P, \L(\underline{l}))\rightarrow
\Pr^k_{G/P}(\L(\underline{l}))(e) \]
is surjective for $1\leq k \leq min\{l_i+1\}$. It follows we get an
exact sequence of locally free
sheaves
\[ 0\rightarrow \Q \rightarrow \H^0(G/P,\L(\underline{l}))\otimes
\O_{G/P} \rightarrow \Pr^k_{G/P}(\L(\underline{l})) \rightarrow 0 \]
on $G/P$. Dualize this sequence to get the sequence
\[  0\rightarrow \Pr^k_{G/P}(\L(\underline{l}))^* \rightarrow
\H^0(G/P,\L(\underline{l}))^* \otimes
\O_{G/P} \rightarrow \Q^* \rightarrow 0 \]
We get a closed immersion
\[ \p(\Q^*)\subseteq \p(\H^0(G/P,\L(\underline{l}))^*) \times G/P \]
of schemes.
Since $\Q$ is locally free  $\Q^*=Coker((T^k)^*)$, hence
by Corollary \ref{irreducible} it follows $D^k(\L(\underline{l}))$
is irreducible and the Theorem is proved.
\end{proof}

The total complex looks as follows:
\[ Tot(\C(T^k)^{i,j})_n=\oplus_{i+j=n}\O(-j)\otimes \pi^*
\H^i(G/P,\wedge^j \Pr^k(\L(\underline{l})^*) .\]
Since $\H^i(G/P,\wedge^j \Pr^k(\L(\underline{l})^*)$ is a finite
dimensional $G$-module there is a decomposition
\begin{equation} \label{decomp}
 \H^i(G/P,\wedge^j \Pr^k(\L(\underline{l})^*) \cong \oplus_\lambda
W_\lambda
\end{equation}
into irreducible $G$-modules. Hence to check if the total complex
$Tot(\C(T^k))^{\bullet, \bullet}$ can be used to give a resolution of
the discriminant
$D^k(\L(\underline{l}))$ one has to calculate the
decomposition \ref{decomp}.


In a series of papers the structure of the jet bundle $\Pr^k(\L)$ of a
line bundle $\L\in \Pic^{\SL(V)}(\SL(V)/P) $ as abstract locally free
sheaf, as left and right $\O$-module and as left and right $P$-module has been
studied (see \cite{maa1},\cite{maa11}, \cite{maa12},
\cite{maa2},\cite{maa3} and \cite{maa31}) using algebraic
techniques, geometric techniques, algebraic group techniques and techniques from universal enveloping
algebras of semi simple Lie algebras. There is work
in progress using techniques similar to the ones introduced in this
paper and the papers \cite{maa1}, \cite{maa2} and \cite{maa3}
on the problem of describing resolutions of ideal sheaves of
$D^k(\L)$ where $\L$ is a line bundle on $\SL(V)/P$ for some parabolic
subgroup $P\subseteq \SL(V)$. This problem may be studied using the
total complex and determinantal schemes.

\end{document}